\documentclass[11pt]{article}
\usepackage[margin=1in]{geometry}
\usepackage{amsmath,amssymb,amsthm,mathtools}
\DeclareMathOperator{\sech}{sech}
\numberwithin{equation}{section}

\usepackage{graphicx}
\usepackage{booktabs}
\usepackage{enumitem}
\usepackage{mathtools}
\usepackage{graphicx}
\usepackage{float}
\usepackage{comment}
\usepackage{tabularx}
\usepackage{hyperref}
\usepackage{booktabs}
\usepackage{array}
\newcolumntype{L}[1]{>{\raggedright\arraybackslash}p{#1}}
\newcolumntype{Y}{>{\raggedright\arraybackslash}X}

\title{\Large\bfseries  On a Mathematical Model Describing Chemotherapeutic\\ Drug Treatment for Tumor Cells}
\author{
\normalsize Xiaoqin Liu\footnote{email: xiaoqin.liu2@wsu.edu} and Hong-Ming Yin\footnote{email: hyin@wsu.edu; Corresponding author}.\\
\small Department of Mathematics and Statistics\\
\small Washington State University, Pullman, WA 99164, USA
}
\date{}
\begin{document}
\maketitle

\begin{abstract}
In this paper, we study a semilinear parabolic PDE system which describes the interaction of normal cells, tumor cells, immune cells, with a chemotherapeutic drug. The model extends the previous model with incorporating strong Allee affects in the normal-tissue and tumor dynamics. Under mild assumptions, we establish global-in-time existence and uniqueness of nonnegative weak solutions and derive $L^{\infty}(Q_{T^*})$ bounds for all $T^*>0$. We then investigate spatiotemporal dynamics of the model and therapy scheduling using an implicit Crank Nicolson Backward Euler (CNBE) scheme. Simulations in a heterogeneous two-dimensional space-dimensional tissue region with three tumor peaks indicate rapid tumor invasion without treatment and significant tumor suppression under pulsed chemotherapeutic treatment. Moreover, in a fixed total dose delivered within the treatment cycle, while keeping each injection duration fixed, concentrated pulses produce stronger early knock-down of tumor density, while more frequent but gentler pulses achieve comparable control of the tumor invasive front while better preserving normal tissue over a four-week period.
\end{abstract}


\section{Introduction}

Cancer remains one of the most severe and complex diseases worldwide, driven by
different biological mechanisms that include genetic mutations, cellular proliferation and
competition, tumor-immune interactions, and tumor response to chemotherapeutic treatment.
Researchers in biology, medicine, engineering, and the mathematical sciences have long sought valuable frameworks to clarify dominant mechanisms, help interpret data, and support therapy design. In particular, mathematical models of cancer growth span a wide range of scales, from non-spatial ordinary differential equations (ODEs) to spatial partial differential equations (PDEs). Previous mathematical modeling works demonstrate that tumor growth, immune interactions, and therapy effects can be represented by dynamical systems that align with clinical observations such as threshold phenomena, tumor dormancy, treatment trade-offs, and oscillations in tumor burden in tumor-immune models arising from delayed immune activation or delayed treatment response (often known as the ``Jeff's phenomenon'' )
(see \cite{Alqahtani2025, MTY2024, CantrellCosner2004, dePillisRadunskaya2001, Eisen1979, EladdadiKimMallet2014,KirschnerPanetta1998,Kuznetsov1994, Kuznetsov1992, Rosen1972, WodarzKomarova2014}).

Our work is rooted in the classical Kuznetsov-Makalkin-Taylor-Perelson ODE model for tumor-immune interactions \cite{Kuznetsov1994}, together with later extensions that include immunotherapy and chemotherapy mechanisms \cite{dePillisRadunskaya2001, KirschnerPanetta1998}. While ODE models capture temporal dynamics, they neglect spatial features that are also important in tumor modeling, such as cell motility, drug penetration. PDE models on the other hand address these limitations by incorporating diffusion and (when appropriate) advection, together with suitable boundary conditions that confine the tissue domain 
(\cite{Andasari2011,ChenLowengrub2011,Ansarizadeh2017,Friedman2005,RooseChapmanPhilip2007,Yin2026}).

The model in this paper is an extension of a four-species reaction-diffusion system including normal tissue density $N(x,t)$, tumor density $T(x,t)$, immune effector density $I(x,t)$, and chemotherapeutic drug concentration $U(x,t)$ introduced by Ansarizadeh, Singh, and Richards \cite{Ansarizadeh2017} in one space-dimension, and later formulated in a higher space-dimension by Yin \cite{Yin2026}. We retain the core biological structure and introduce two modeling modifications to distinguish this work:

\begin{itemize}
\item \textbf{Strong Allee effects in $N$ and $T$.}
We incorporate Allee-type growth in both the normal and tumor equations, so that each population may exhibit negative net growth below a critical threshold value. The term ``Allee effect'' was originally used to identify biological mechanisms that reduce per-capita growth at low density in population dynamics (see \cite{Courchamp2008,DrakeKramer2011,Stephens1999}), and later, it has also been adopted in cellular modeling \cite{Bottger2015,DelitalaFerraro2020,Gerlee2022,Johnson2019}. In particular, a strong Allee effect means that there exists a threshold density below which the net growth becomes negative and extinction is inevitable, whereas densities above the threshold can grow and persist.

\item \textbf{General second-order operators in divergence form.} To allow more flexible spatial dynamics, we formulate each equation with a general second-order elliptic operator in divergence form (with a first-order advection term). Biologically, advection is particularly important for the drug variable, since drug delivery from blood vessels to tumor tissue is not purely diffusive and is, in fact, mainly affected by directed delivery mechanisms \cite{DewhirstSecomb2017}. For the sake of analysis, we allow advection terms in all four equations so that the well-posedness theory applies to the full operator family. However, in the numerical simulations, we set the advection terms in the $N$, $T$, and $I$ equations to zero, while having a nonzero advection field only in the $U$ equation, reflecting diffusion-dominated cell motility and advection-driven drug transport.
\end{itemize}

\subsection{Final PDE model}
Motivated by the above biological facts and the four component model in \cite{Yin2026}, let $\Omega\subset\mathbb{R}^n$ be a bounded domain with smooth boundary
$\partial\Omega$, and set $Q_{T^*}:=\Omega\times(0,T^*)$, and $S_{T^*}:=\partial\Omega\times(0,T^*)$, where $T^*>0$ denotes time. We also denote $(u_1,u_2,u_3,u_4)=(N,T,I,U)$, representing normal tissue density, tumor density,
immune effector cells density, and chemotherapeutic drug concentration, respectively. For $k=1,\dots,4$, define the divergence-form operator
\begin{equation}\label{lk}
L_k[u_k]
:= -\sum_{i,j=1}^n \partial_ {x_j}(d_k^{ij}(x,t)\,
\partial_{x_i} u_k)
+ \sum_{i=1}^n h_k^{i}(x,t)\,\partial_{x_i} u_k,
\end{equation}
where the coefficients satisfy
\begin{equation}\label{lkcoeff}
d_k^{ij}(x,t),\; h_k^{i}(x,t)\in L^\infty\!\big(Q_{T^*}\big),
\qquad \text{for } k=1,\dots,4.
\end{equation}
We have reached our extended model:
\begin{align}
N_t + L_1[N] &= F_1(N,T,I,U), && (x,t)\in Q_{T^*}, \label{sys_1}\\
T_t + L_2[T] &= F_2(N,T,I,U), && (x,t)\in Q_{T^*}, \label{sys_2}\\
I_t + L_3[I] &= F_3(N,T,I,U), && (x,t)\in Q_{T^*}, \label{sys_3}\\
U_t + L_4[U] &= F_4(N,T,I,U), && (x,t)\in Q_{T^*}, \label{sys_4}
\end{align}
subject to homogeneous Neumann boundary conditions
\begin{equation}\label{sys_b}
(D_1(x,t)\nabla N)\cdot\boldsymbol{\nu}
=(D_2(x,t)\nabla T)\cdot\boldsymbol{\nu}
=(D_3(x,t)\nabla I)\cdot\boldsymbol{\nu}
=(D_4(x,t)\nabla U)\cdot\boldsymbol{\nu}
=0,
\end{equation} for all $(x,t) \in S_{T^*}$, 
and initial conditions
\begin{equation}\nonumber
N(x,0)=N_0(x)=u_{01}(x),\quad
T(x,0)=T_0(x)=u_{02}(x),
\end{equation}
\begin{equation}\label{sys_i}
I(x,0)=I_0(x)=u_{03}(x),\quad
U(x,0)=U_0(x)=u_{04}(x),
\end{equation}for all $x\in\Omega$. Here $D_k(x,t)$ are diffusion matrices whose $i$-th row, $j$-th column entry is denoted by $d_k^{i,j}(x,t)$ for each $k$, and 
$$
D_k(x,t)\nabla u_k\cdot\boldsymbol{\nu} = \sum_{i,j=1}^n d_k^{ij}(x,t)\,\partial_{x_i}u_k\,\nu_j,
$$ where $\boldsymbol{\nu}$ is the outward unit
normal vector on $\partial\Omega$. The nonlinear reaction terms are
\begin{align}
F_1(N,T,I,U)
&= r_2 N\big(1-b_2N\big)\Big(\tfrac{N}{A_2}-1\Big)
- c_4TN - a_3\big(1-e^{-U}\big)N, \label{ode_1}\\
F_2(N,T,I,U)
&= r_1 T\big(1-b_1T\big)\Big(\tfrac{T}{A_1}-1\Big)
- c_2IT - c_3TN - a_2\big(1-e^{-U}\big)T, \label{ode_2}\\
F_3(N,T,I,U)
&= s(x,t)+\frac{\rho IT}{\alpha+T}
- c_1IT - k_1I - a_1\big(1-e^{-U}\big)I, \label{ode_3}\\
F_4(N,T,I,U)
&= v(x,t)\,H(N-a_0) - k_2U. \label{ode_4}
\end{align}

In our model, tumor and normal cells follow logistic growth with constant growth rates $r_1,r_2$ and with carrying capacities $1/b_1,1/b_2$, respectively, and according to our modification, this growth for each population experiences negative growth when it falls below its corresponding Allee threshold $A_1$ or $A_2$. Interaction rates are represented by $c_1$ (loss of immune cells due to tumor cells), $c_2$ (loss of tumor cells due to immune cells), $c_3$ (loss of tumor cells due to competition between tumor and normal cells), and $c_4$ (damage of normal cells by tumor cells). Chemotherapeutic drug action enters through fractional-kill terms $a_k(1-e^{-U})$ for $k=1,2,3$. The immune equation includes a baseline influx $s(x,t)$ (tumor-independent) and a tumor-driven proliferation term $\frac{\rho I T}{\alpha+T}$, with natural decay at rate $k_1$. Drug dynamics are driven by $v(x,t) H(N-a_0)$ (delivery only when $N>a_0$ for some positive safety threshold constant $a_0$), with linear decay at rate $k_2$.

This work has two main parts. In Section $2$, under mild structural assumptions, we establish global well-posedness of weak solutions. We first obtain a unique local-in-time weak solution via the Schauder fixed-point theorem and then, using a blow-up alternative, extend it to a global-in-time solution once a priori $L^\infty(Q_{T^*})$ bounds for $(u_1,\cdots,u_4)$ are obtained. In Section $3$, we develop a CNBE scheme. We then perform numerical simulations on a heterogeneous two space-dimensional tissue with three tumor peaks to illustrate spatial tumor progression and to investigate pulsed chemotherapy drug delivery design. In particular, for fixed total dosage drug treatment (i.e. $V_0\tau N^*$ is held fixed), we vary the injection rate $V_0$, injection duration $\tau$, and number of injections $N^*$ to demonstrate rate-duration-frequency trade-offs in treatment scheduling.

\section{Global well-posedness}
In this section, we establish global well-posedness of weak solutions to the system \eqref{sys_1}-\eqref{sys_4} under some assumptions that also consistent with the biological setting. In order to state the main result precisely, we first introduce the notation and give the definition of a weak solution for our model. We then list the hypotheses used throughout the paper.

\subsection{Preliminary}
The notations and function-space conventions used in this work follow those of Evans \cite{Evans2010} and Yin \cite{YinPDEBook}. Throughout the paper, we use $C$, $M$, etc. to denote generic positive constants depending only on known data (their values may change from line to line). We use $C_i$, $M_i$, etc. (or $C_{i,j}$, $M_{i,j}$, etc.) to denote distinct indexed fixed constants. Similarly, $\epsilon>0$ denotes a generic small parameter, and $\epsilon_i$ or $\epsilon_{i,j}$ to indicate independent small parameters chosen at different steps. When the dependence of a constant matters, we write $C=C(\cdot)$ or $M=M(\cdot)$. We use boldface to indicate vectors, and scalar values appear in plain type.

For $0<\beta\le 1$, the parabolic H\"older space $C^{\beta,\beta/2}(\overline{Q_{T^*}})$ is defined by
$$
C^{\beta,\beta/2}(\overline{Q_{T^*}})
:=\Big\{u:\overline{Q_{T^*}}\to\mathbb{R}\ \big|\ 
\|u\|_{C^{\beta,\beta/2}(\overline{Q_{T^*}})}<\infty\Big\},
$$
with norm
$$
\|u\|_{C^{\beta,\beta/2}(\overline{Q_{T^*}})}
:= \|u\|_{L^\infty(\overline{Q_{T^*}})}
+ \sup_{\substack{(x_1,t_1)\neq (x_2,t_2)\\ (x_i,t_i)\in\overline{Q_{T^*}}}}
\frac{|u(x_1,t_1)-u(x_2,t_2)|}{|x_1-x_2|^\beta + |t_1-t_2|^{\beta/2}}.
$$
Let $X$ be a Banach space and let $1\le p<\infty$. The Bochner space $L^p((0,T^*);X)$ is defined as:
$$
L^p((0,T^*);X)
:= \Big\{u(x,t): \forall t \in [0,T^*], u(x,t) \in X\Big\},
$$
and it is equipped with the norm
$$
\|u\|_{L^p((0,T^*);X)}
:= \Big(\int_0^{T^*}\|u(t)\|_X^p\,dt\Big)^{1/p}.
$$
Define
$$
V_2(Q_{T^*})
:= L^2\big((0,T^*);H^1(\Omega)\big)\cap L^\infty\big((0,T^*);L^2(\Omega)\big),
$$
equipped with the norm
$$
\|u\|_{V_2(Q_{T^*})}
:= ess\sup_{0<t<T^*}\|u(t)\|_{L^2(\Omega)}
+ \|\nabla u\|_{L^2(Q_{T^*})},
$$ where
$$
\|\nabla u\|_{L^2(Q_{T^*})}
:= \Big(\int_0^{T^*}\!\!\int_\Omega |\nabla u(x,t)|^2\,dx\,dt\Big)^{1/2}.
$$

For convenience, the spaces $C(Q_{T^*};\mathbb{R}^m)$, $L^\infty(Q_{T^*};\mathbb{R}^m)$, etc.\
consist of vector-valued functions $\mathbf{u}=(u_1,\dots,u_m):Q_{T^*}\to\mathbb{R}^m$ such that
$u_k\in C(Q_{T^*})$, $L^\infty(Q_{T^*})$, etc.\ for $k=1,\dots,m$.

\medskip
\noindent\textbf{Definition 2.1.1} (Weak solution for Neumann boundary condition).\\
We say that $\mathbf{u}$ is a weak solution to \eqref{sys_1}–\eqref{sys_4} on $Q_{T^*}=\Omega\times(0,T^*)$ if, for each $k=1,\dots,4$,
$$
u_k \in V_2(Q_{T^*}) \,\cap\, L^\infty(Q_{T^*}),
$$
and for every test function
$$
\phi_k \in L^2\!\big(0,T^*; H^1(\Omega)\big)\,\cap\, H^1\!\big(0,T^*; L^2(\Omega)\big)
\text{ with } \phi_k(\cdot,T^*)=0,
$$
the integral identity
\begin{equation}\label{definition}
\begin{aligned}
&\iint_{Q_{T^*}} \big(-u_k\,\partial_t \phi_k\big)\,dx\,dt
\,+\, \iint_{Q_{T^*}} \sum_{i,j=1}^{n} d_k^{ij}(x,t)\,\partial_{x_i} u_k\,\partial_{x_j}\phi_k \,dx\,dt+\, \iint_{Q_{T^*}} \sum_{i=1}^{n} h_k^{\,i}(x,t)\,\partial_{x_i} u_k\,\phi_k \,dx\,dt\\ &= \iint_{Q_{T^*}} F_k\!\big(x,t,u_1,\dots,u_4\big)\,\phi_k \,dx\,dt
\,+\, \int_{\Omega} u_{0k}(x)\,\phi_k(x,0)\,dx
\end{aligned}
\end{equation}
holds. We call $\mathbf{u}$ a local-in-time weak solution if \eqref{definition} holds on $Q_{T^*}$ for some $T^*>0$, and a global-in-time weak solution if it holds for any $T^*>0$.

\subsection{Hypothesis}
For this work, we impose the following hypotheses:

\begin{comment}
In \textbf{H(1)} we assume that the divergence-form elliptic operators are uniformly elliptic with bounded coefficients, an assumption appropriate for the parabolic system studied here. \textbf{H(2)} enforces biological sign constraints: all logistic growth rates, interaction coefficients, and Allee thresholds are positive. \textbf{H(3)} requires the immune influx $s(x,t)$ and drug input capacity $v(x,t)$ to be bounded and nonnegative. \textbf{H(4)} employs a $C^1$ smoothed Heaviside-like function to model on–off gating while preserving differentiability for analysis. Finally, \textbf{H(5)} prescribes nonnegative, Hölder-regular initial data, reflecting physically meaningful, bounded spatial profiles at $t=0$ and ensuring compatibility with the well-posedness framework.
\end{comment}

\begin{itemize}
\item\textbf{H(1)}: For $k=1,2,3,4$, $D_k(x,t)=(d_k^{ij}(x,t))_{n\times n}$ is a positive definite matrix satisfying 
$$
d_0|\xi|^2 \leqslant \sum^n_{i,j=1}d_k^{ij}(x,t)\xi_i\xi_j \leqslant d_1|\xi|^2,\quad \text{ for all } \xi \in \mathbb{R}^n,
$$ 
for some constants $d_0, d_1>0$. Moreover, there exist $h^*$ such that, for $k=1,2,3,4$,
$$
\|h^{i}_k(x,t)\|_{L^{\infty}(Q_{T^*})}\leqslant h^*, \text{ for each }i=1,\cdots n.
$$
\item\textbf{H(2)}: All parameters $a_0, \cdots, a_3, b_1, b_2, c_1, \cdots, c_4, r_1, r_2, k_1,k_2, A_1, A_2$ and $\rho, \alpha$ are positive constants.

\item\textbf{H(3)}: $s(x,t)$ and $v(x,t)$ are nonnegative over $Q_{T^*}$ and belong to $L^\infty(Q_{T^*})$ for all $T^*>0$.

\item\textbf{H(4)}: The Heaviside-like function in \eqref{ode_4} satisfies $H(y) \in C^1(\mathbb{R})$, and $H(y)=0$ in $(-\infty,0]$, $H(y)=1$ in $[\delta,\infty)$ for any fixed small $\delta>0$.

\item\textbf{H(5)}: The initial condition for each species 
$$
\begin{aligned}
&N_0(x),\,T_0(x),\,I_0(x),\,U_0(x) \in C^{\alpha}(\overline{\Omega}), \text{ and }\\
&N_0(x),\,T_0(x),\,I_0(x),\,U_0(x) \geqslant 0 \quad \text{for all } x \in \Omega .
\end{aligned}
$$
\end{itemize}

\subsection{Main results}

\medskip
\noindent\textbf{Theorem 2.3.1} (Global existence and uniqueness of the weak solution) Suppose \textbf{H(1)}–\textbf{H(5)} hold. Then for any $T^*>0$, the system \eqref{sys_1}-\eqref{sys_4} admits an unique global-in-time weak solution $\mathbf{u}=(u_1,u_2,u_3,u_4)$ on $Q_{T^*}$ in the sense of Definition 2.1.1, with  
$$ 
\mathbf{u} \in V_2(Q_{T^*}; \mathbb{R}^4) \cap L^{\infty}(Q_{T^*}; \mathbb{R}^4) \cap C^{\beta, \frac{\beta}{2}}(\overline{Q_{T^*}}; \mathbb{R}^4).
$$ Moreover, there exists a constant $M(T^*)$ such that 
$$
0 \le u_1(x,t),u_2(x,t),u_3(x,t),u_4(x,t)\le M(T^*)
$$ over $(x,t) \in Q_{T^*}$.

\medskip
We now proceed to derive some useful a priori estimates for the system \eqref{sys_1}-\eqref{sys_4}. In particular, to obtain an $L^\infty$ estimate for the immune density
$I(x,t)$, we will invoke the $L^1(Q)$-$L^\infty(Q)$ estimate established by Yin \cite{Yin2026}; we briefly indicate below how its hypotheses are verified in our setting. In Appendix A, we also give a slightly different proof to \textbf{(b)} of Lemma 2.3.2.

\medskip
\noindent\textbf{Lemma 2.3.2} ($L^2(Q_{T^*})$ to $L^\infty(Q_{T^*})$ estimate).
Let $u$ be a weak solution of the below scalar initial value boundary problem:
\begin{align} 
&\partial_t u + L[u] = f(x,t,u), && (x,t)\in Q_{T^*}, \label{pdes1}\\ 
&(D(x,t)\nabla u)\cdot\nu = 0, && (x,t)\in S_{T^*},\label{pdes2}\\ 
&u(x,0) = g(x), && x\in\Omega,\label{pdes3} \end{align} 
where $L[\cdot]$ and the diffusion matrix $D(x,t)$ satisfy \textbf{H(1)} with $k=1$. Assume $g(x) \ge 0$ in $\Omega$, and $f(x,t,0)\ge 0$ and the one-sided growth condition,
$$
f(x,t,u)\,u \le \bar C\,(1+u^2)\qquad \text{for all } u\ge 0,\ \text{over } (x,t)\in Q_{T^*}.
$$

 \medskip \noindent\textbf{(a) $L^2$–estimate.} 
 For every $T^*>0$ there exists $M_0(T^*)>0$, depending only on $T^*$, $\Omega$, the coefficients of the operator $L[\cdot]$, $\bar C$ and $\|g(x)\|_{L^\infty(\Omega)}$, such that 
 $$ 
 \sup_{0\le t\le T^*}\int_\Omega u(x,t)^2\,dx \le M_0(T^*). 
 $$ 
 
 \medskip \noindent\textbf{(b) $L^\infty$–estimate }
 Suppose \textbf(a) holds, then there exists $\hat{M}(T^*)>0$, depending only on $M_0(T^*)$, $\Omega$, and the coefficients of $L[\cdot]$, such that 
 $$
 \sup_{0\le t\le T^*}\|u(\cdot,t)\|_{L^\infty(\Omega)} \le \hat{M}(T^*). 
 $$

\medskip
\noindent\textit{Proof.}
First of all, since $g(x)\ge 0$ and $f(x,t,0)\ge 0$, maximum principle implies $u(x,t)\ge 0$ on $Q_{T^*}$, so the assumption $f(x,t,u)\,u \le \bar C\,(u^2+1)$ applies.

To show \textbf(a), multiplying \eqref{pdes1} by $u$ and using \textbf{H(1)} and the one-sided growth condition
yields an $L^2$ estimate
$$
\sup_{0\le t\le T^*}\|u(\cdot,t)\|_{L^2(\Omega)}^2 \le M_0(T^*),
$$
for some $M_0(T^*)>0$ depending only on known data and $T^*$.

To prove \textbf(b), since $\Omega$ is bounded, Hölder's inequality gives, for any $t\in[0,T^*]$,
$$
\int_\Omega u(x,t)\,dx \le |\Omega|^{1/2}\|u(\cdot,t)\|_{L^2(\Omega)}
\le |\Omega|^{1/2} M_0(T^*)^{1/2}.
$$
Hence 
$$
\sup_{0\le t\le T^*}\int_\Omega u(x,t)\,dx<\infty.
$$ Therefore the hypotheses of
Lemma 2.1 in \cite{Yin2026} are satisfied (restricted to $Q_{T^*}$), and we obtain
$$
\sup_{0\le t\le T^*}\|u(\cdot,t)\|_{L^\infty(\Omega)}\le \hat M(T^*),
$$ the desired $L^\infty$ estimate on $Q_{T^*}$ depends only on known data and $T^*$. \hfill$\square$

\subsection{A priori estimates}

\noindent\textbf{Lemma 2.4.1} (A priori $L^{\infty}(Q_{T^*})$ bounds for N, T, I, U for any $T^*>0$). Assume that hypotheses \textbf{H(1)}–\textbf{H(5)} hold for the system \eqref{sys_1}-\eqref{sys_4} with initial data \eqref{sys_i}.
Then for any $T^*>0$, there exists a constant $M(T^*)>0$, depending only on known data, such that 
$$
0 \,\le\, N(x,t),\, T(x,t),\, I(x,t), \,U(x,t) \,\le\, M(T^*) 
\quad \text{for all } (x,t)\in Q_{T^*}.
$$

\noindent\textit{Proof.}
We first show the nonnegativity of $N,T,I,U$. By \textbf{H(3)}–\textbf{H(4)},
\begin{equation*}
\begin{split}
&F_1(x,t,0,T,I,U)=0,\\ 
&F_2(x,t,N,0,I,U)=0, \\
&F_3(x,t,N,T,0,U)=s(x,t),\\
&F_4(x,t,N,T,I,0)= v(x,t)\,H(N-a_0) \,\ge\, 0,
\end{split}
\end{equation*}
and by \textbf{H(5)} the initial data satisfy $N_0,\,T_0,\,I_0,\,U_0 \ge 0$. By applying the parabolic maximum principle under homogeneous Neumann boundary conditions and a componentwise comparison principle, we see $\underline{\mathbf{u}}\equiv \mathbf{0}$ is a lower solution for $(N,T,I,U)$. Hence
\begin{equation}\label{sys_lower}
N(x,t),\,T(x,t),\,I(x,t),\,U(x,t) \,\ge\, 0 \quad \text{for all } (x,t)\in Q_{T^*}.
\end{equation}

\emph{Upper bound for $N$.}
From \eqref{ode_1} and the nonnegativity of $T$, $N$, and $(1-e^{-U})$ and \textbf{H(2)},
\begin{equation}\label{upper_N}
\begin{aligned}
F_1(N,T,I,U)
&= r_2 N(1-b_2 N)\Big(\frac{N}{A_2}-1\Big) - c_4 T N - a_3(1-e^{-U})N\\
&\le\, r_2 N(1-b_2 N)\Big(\frac{N}{A_2}-1\Big).
\end{aligned}
\end{equation}
Choose
$$
M_1 \,\ge\, \max\Big\{\tfrac{1}{b_2},\, A_2,\, \|N_0\|_{L^\infty(\Omega)}\Big\}
$$
so that $F_1(M_1,T,I,U)\le 0$, therefore, $M_1$ is an upper solution such that
\begin{equation} \label{M_1}
 0\le\,N(x,t) \,\le\, M_1 \quad\text{on }Q_{T^*}.
\end{equation}

\emph{Upper bound for $T$.}
Similarly, from \eqref{ode_2},
\begin{equation}\label{upper_T}
\begin{aligned}
F_2(N,T,I,U)
&\le\, r_1 T(1-b_1 T)\Big(\frac{T}{A_1}-1\Big).
\end{aligned}
\end{equation}
Choose
$$
M_2 \,\ge\, \max\Big\{\tfrac{1}{b_1},\, A_1,\, \|T_0\|_{L^\infty(\Omega)}\Big\}
$$
so that $F_2(N,M_2,I,U)\le 0$, hence $M_2$ is an upper solution and 
\begin{equation} \label{M_2}
 0\le\,T(x,t) \,\le\, M_2 \quad\text{on }Q_{T^*}.
\end{equation}

\emph{Upper bound for $U$.}
From \eqref{ode_4},
\begin{equation}\label{upper_U}
\begin{aligned}
F_4(N,T,I,U) &= v(x,t)H(N-a_0)-k_2 U \,\le\, \|v\|_{L^\infty(Q_{T^*})} - k_2 U.
\end{aligned}
\end{equation}
Choose
$$
M_3 \,\ge\, \max\Big\{\|U_0\|_{L^\infty(\Omega)},\, \tfrac{1}{k_2}\|v\|_{L^\infty(Q_{T^*})}\Big\}
$$
so that $F_4(N,T,I,M_3)\le 0$, and thus $M_3$ is an upper solution and
\begin{equation} \label{M_3}
0\le\,U(x,t) \,\le\, M_3 \quad\text{on }Q_{T^*}.
\end{equation}

\emph{Upper bound for $I$.}
From \eqref{sys_3} and \eqref{ode_3}, $I$ satisfies
\begin{equation}\label{I_eq}
\partial_t I + L_3[I]
= s(x,t) + \frac{\rho IT}{\alpha+T} - k_1 I - c_1 IT - a_1(1-e^{-U})I,
\end{equation}
where
$$
L_3[I]=-\sum_{i,j=1}^n \partial_{x_j}\!\big(d_3^{ij}(x,t)\partial_{x_i}I\big)
+\sum_{i=1}^n h_3^i(x,t)\partial_{x_i}I.
$$

Multiplying \eqref{I_eq} by $I$ and using \textbf{H(1)}, integration by parts yields
\begin{equation}\label{I_energy0}
\frac12\frac{d}{dt}\int_\Omega I^2\,dx
+d_0\int_\Omega |\nabla I|^2\,dx
+k_1\int_\Omega I^2\,dx
\le \underbrace{h^*\int_\Omega |\nabla I|\,I\,dx}_{I_{1}}
+\underbrace{\int_\Omega s(x,t)\,I\,dx}_{I_2}
+\rho\int_\Omega I^2\,dx.
\end{equation}
By Young's inequality (choose $\epsilon_1=\tfrac{d_0}{2h^*}$),
$$
I_{1}
\le \frac{d_0}{2}\int_\Omega |\nabla I|^2\,dx
+\frac{(h^*)^2}{2d_0}\int_\Omega I^2\,dx,
$$
By Young's inequality with $\epsilon$
$$
I_2
\le \epsilon \int_\Omega I^2\,dx + C(\epsilon)\,|\Omega|.
$$
Substituting into \eqref{I_energy0} gives
\begin{equation}\label{I_energy1}
\frac{d}{dt}\int_\Omega I^2\,dx +C_1\int_\Omega I^2\,dx
\le C_2,
\end{equation}
for constants $C_1 \,=\, k_1 - \rho - \frac{(h^*)^2}{2d_0}$ and $C_2>0$ depending only on $d_0,h^*,\rho,k_1,\|s\|_{L^\infty(Q_{T^*})},|\Omega|$. Therefore, by Gr\"onwall's inequality, for any $T^*>0$,
\begin{equation}\label{I_L2_bound}
\sup_{0\le t\le T^*}\int_\Omega I(x,t)^2\,dx \le \hat M_4(T^*).
\end{equation}

Moreover, we verify the one sided growth condition for $F_3(N,T,I,U)$. Since \textbf{H(2)} and the nonnegativity of $I$, $T$, and $U$,
\begin{equation}\label{F_I}
\begin{aligned}
F_3(N,T,I,U)\,I
&= \Big(s(x,t) + \frac{\rho IT}{\alpha+T} - c_1 IT - k_1 I - a_1(1-e^{-U})I\Big)I\\
&\le \Big(s(x,t) + \frac{\rho IT}{\alpha+T}\Big)I.
\end{aligned}
\end{equation}
Using $0 \le T \le M_2$ and so $0 \le \tfrac{T}{\alpha+T} \le 1$, we obtain a constant $\bar C>0$ depending only on $\|s\|_{L^\infty(Q_{T^*})}$ and $\rho$ such that
$$
F_3(N,T,I,U)\,I
\le \|s\|_{L^\infty(Q_{T^*})}I + \rho I^2
\le \bar C\,(I^2+1),
$$ the last inequality is obtained through Young's inequality. Therefore, Lemma~2.3.2 implies a constant $M_4(T^*)$ depends only on the known data and $T^*$ such that
$$
\sup_{0 \le t \le T^*}\|I(\cdot,t)\|_{L^{\infty}(\Omega)} \,\le\, M_4(T^*).
$$

Finally, set $M(T^*):=\max\{M_1, M_2, M_3, M_4(T^*)\}$ to conclude the stated bounds for $N,T,I,U$ above. \hfill$\square$

\subsection{Global Existence and Uniqueness of a Semilinear Parabolic PDE System}
\medskip
In order to prove Theorem 2.3.1, we first show the system has unique local-in-time weak solution, we will then extend it to a global one.

\medskip
\noindent\textbf{Lemma 2.5.1} (Local-in-time existence and uniqueness).
Suppose that \textbf{H(1)}-\textbf{H(5)} hold. Then there exists $T^*>0$ sufficiently small such that \eqref{sys_1}-\eqref{sys_4} admits a unique local-in-time weak solution in the sense of Definition 2.1.1, with
$$
\mathbf{u}=(u_1,u_2,u_3,u_4)\in V_2(Q_{T^*};\mathbb{R}^4)\cap L^\infty(Q_{T^*};\mathbb{R}^4)\cap
C^{\beta,\beta/2}(\overline{Q_{T^*}};\mathbb{R}^4)
$$
for some $0<\beta\le 1$.

\medskip
\noindent\textit{Proof.}
Fix $0<T^*\le 1$ (to be chosen) and set $$
K_0:=\sum_{k=1}^4\|u_{0k}\|_{L^\infty(\Omega)}+1,
$$ let 
$$
K:=\Big\{\mathbf{w}=(w_1,\cdots,w_4)\in L^\infty(Q_{T^*};\mathbb{R}^4):\  \|\mathbf{w}\|_{L^{\infty}(Q_{T^*}; \mathbb{R}^4)}\le K_0\Big\}.
$$
For given $\mathbf{w}\in K$, consider the linear system for $k=1,\dots,4$,
\begin{align}
&\partial_t u^{*}_k\,+\, L_k[u^{*}_k] \,=\, F_k(x,t,\mathbf{w}), && (x, t) \in Q_{T^*}, \label{linearpde1}\\
&(D_k(x,t)\,\nabla u^{*}_k)\cdot\nu\,=\,0, && (x, t) \in S_{T^*}, \label{linearpde2}\\
&u^{*}_k(x, 0) \,=\, u_{0k}(x), && x \in \Omega. \label{linearpde3}
\end{align}
By standard linear parabolic theory \cite{Evans2010,Ladyzhenskaya1968}, suppose \textbf{H(1)}-\textbf{H(5)} hold, then there exists a unique weak solution $\mathbf{u}^*\in V_2(Q_{T^*};\mathbb{R}^4)$.
Moreover, by the De Giorgi-Nash estimates and \textbf{H(5)},
\begin{equation}\label{DN_est}
\begin{aligned}
\|\mathbf{u}^*\|_{C^{\beta,\beta/2}(\overline{Q_{T^*}};\mathbb{R}^4)}\le C,
\end{aligned}
\end{equation}
for some constant $C$ depends only on $K_0$ and known data.

Define the map $M:K\to L^\infty(Q_{T^*};\mathbb{R}^4)$ by $M[\mathbf{w}]=\mathbf{u}^*$.
Using \eqref{DN_est} with $(\hat x,\hat t)=(x,0)$, we choose $
T^* = \min\Big\{\Big(\tfrac{1}{2C(T^*)}\Big)^{\!2/\beta},\,1\Big\}$ to obtain
$$
\|\mathbf{u}^*\|_{L^\infty(Q_{T^*})}
\le \sum_{k=1}^4\|g_k\|_{L^\infty(\Omega)} + C (T^*)^{\beta/2} \le K_0.
$$ Therefore $M:K\to K$.

Compactness of $M$ follows because $M[\mathbf{w}]$ is bounded in
$C^{\beta,\beta/2}(\overline{Q_{T^*}};\mathbb{R}^4)$ and the embedding
$C^{\beta,\beta/2}(\overline{Q_{T^*}})\hookrightarrow L^\infty(Q_{T^*})$ is compact
(Arzel\`a-Ascoli Theorem). 

Let $\mathbf{w}^j=(w_1^j,\cdots,w_4^j)$. For each $j\ge 1$ set $\mathbf{u}^{*j}=M[\mathbf{w}^j]$ and $\mathbf{u}^*=M[\mathbf{w}]$. Define $\mathbf{v}^j\,:=\,(v^j_1,\cdots,v^j_4)$, where $v^j_k= u_k^{*j}-u_k^*$ for $k=1,\cdots,4$. Continuity of $M$ follows from an $L^2$-energy estimate for the difference $\mathbf{v}^j$ together with Gr\"onwall's inequality, and the fact that $F_k(x,t,\mathbf{w})$ are locally Lipschitz in $\mathbf{w}$ for all $k=1, \cdots,4$. A detailed prove of this part can be found in Appendix A.

Therefore, by the Schauder fixed-point theorem, the mapping $M$ admits a fixed point $\hat{\mathbf{u}}\in K$, which is a local-in-time weak solution to \eqref{sys_1}-\eqref{sys_4}.
Uniqueness follows by applying a standard $L^2$-energy estimate to the difference of two weak solutions $u_k'$ and using Gr\"onwall's inequality, again invoking local Lipschitz continuity of $F_k$, one can show $u_k'=0$ for all $k=1,\cdots,4$, proving uniqueness. \hfill$\square$

\medskip
We now extend the local-in-time weak solution obtained in Lemma 2.5.1 to a unique global-in-time weak solution via a blow-up alternative. To state it precisely, we first define the maximal existence time of the local-in-time weak solution by
\begin{equation}\label{Tmax}
\begin{aligned}
T^*_{\max}
&:= \sup\Big\{\,T^*>0 \ \Big|\ 
\text{there exists a unique local-in-time weak solution } \mathbf{u} \\
&\qquad\text{to \eqref{sys_1} -\eqref{sys_4} on }\Omega\times(0,T^*)\Big\}.
\end{aligned}
\end{equation}

\medskip
\noindent\textbf{Corollary 2.5.2 (Blow-up alternative).}
Assume \textbf{H(1)}-\textbf{H(5)}. Let $\mathbf{u}$ be the unique local-in-time weak solution on $Q_{T^*}$ for some $T^*>0$, and let $T^*_{\max}$ be defined as above. Then exactly one of the following holds:

\medskip
\noindent\textbf{(B1)} $T^*_{\max}<\infty$ and the solution cannot exist beyond $T^*_{\max}$, moreover,
$$
\limsup_{t\rightarrow T^{*-}_{\max}} \sum_{k=1}^4 \|u_k(\cdot,t)\|_{L^\infty(\Omega)} = \infty.
$$

\noindent\textbf{(B2)} The solution can be extended to $Q_{T^*}$ for any $T^*>0$.

\medskip
\noindent\textit{Proof.}
Assume $T^*_{\max}<\infty$ and
$$
M := \limsup_{t\rightarrow T^{*-}_{\max}} \sum_{k=1}^4 \|u_k(\cdot,t)\|_{L^\infty(\Omega)} < \infty.
$$
We choose an increasing sequence of times $\{t_i\}_{i\in\mathbb{N}}$ with $t_i \rightarrow T^{*-}_{\max}$ such that
$$
\sum_{k=1}^4 \|u_k(\cdot,t_i)\|_{L^\infty(\Omega)} \le C_1
$$
for some constant $C_1>0$. By Lemma 2.5.1, there exist $\beta\in(0,1]$ and $C_2>0$ (independent of $i$) such that
$$
\|u_k(\cdot,t_i)\|_{C^\beta(\overline{\Omega})}\le C_2,
$$ for all $k=1,\dots,4$.
Hence for each $k$, the sequence $\{u_k(\cdot,t_i)\}$ is uniformly bounded and equicontinuous in $C^\beta(\overline{\Omega})$. By the Arzelà-Ascoli theorem, we can extract a subsequence (using a diagonal argument, so that this single subsequence guarantees convergence for all $k$), still denoted by $\{t_i\}$, such that for each $k=1,\cdots,4$,
$$
\{u_k(\cdot, t_i)\} \rightarrow u^*_k \quad \text{uniformly on } \overline{\Omega},
$$
with $u^*_k \in C^{\beta}(\overline{\Omega})$. On the other hand, by the Hölder continuity in time from Lemma 2.5.1, we have that the limit
$$
u_k(\cdot,T^{*}_{\max}) := \lim_{t\rightarrow T^{*-}_{\max}} u_k(\cdot,t)
$$
exists in $C^{\beta}(\overline{\Omega})$. By uniqueness of the limit, we must have $u^*_k = u_k(\cdot,T^{*}_{\max})$.  
Therefore,
$$
\{u_k(\cdot, t_{i})\} \rightarrow u_k(\cdot,T^{*}_{\max}) \quad \text{uniformly on } \overline{\Omega}.
$$
In particular,
$$
\sum_{k=1}^4 \|u_k(\cdot, T^{*}_{\max})\|_{L^{\infty}(\Omega)} \le C_2.
$$

Now take $\mathbf{u}(x,T^{*}_{\max})=(u_1(x,T^{*}_{\max}),\cdots,u_4(x,T^{*}_{\max}))$ as a new initial value. By Lemma 2.5.1, there exists $\delta>0$ such that the system \eqref{sys_1}-\eqref{sys_4} admits a weak solution on $[T^{*}_{\max},T^{*}_{\max}+\delta)$. This contradicts the maximality of $T^*_{\max}$. Therefore, if $T^*_{\max}<\infty$, the $L^\infty$-norm must blow up as $t\rightarrow T^*_{\max}$, proving \textbf{(B1)}. Otherwise \textbf{(B2)} holds. \hfill$\square$

\subsection{Proof of Theorem~2.3.1}
\noindent\textit{Proof.} 
\begin{comment}
We note that the reaction terms $F_k(x,t,\mathbf{u})$ in \eqref{ode_1}-\eqref{ode_4} are locally Lipschitz in
$\mathbf{u}=(N,T,I,U)$ for all $k=1,\cdots,4$. Therefore the local existence and uniqueness result
of Lemma 2.5.1 applies to \eqref{sys_1}-\eqref{sys_4}.
\end{comment}
By Lemma 2.5.1, there exists $T_0>0$ such that \eqref{sys_1}-\eqref{sys_4} admits a unique weak solution
$\mathbf{u}$ on $Q_{T_0}$. Let $T^*_{\max}$ be the maximal existence time defined in \eqref{Tmax}.
By the blow-up alternative Corollary 2.5.2, if $T^*_{\max}<\infty$ then it must be true that
$$
\limsup_{t\to T^*_{\max}-}\sum_{k=1}^4 \|u_k(\cdot,t)\|_{L^\infty(\Omega)}=\infty.
$$
On the other hand, Lemma 2.4.1 gives the a priori estimates
$$
\|u_k\|_{L^\infty(Q_{T^*})}\le M(T^*)<\infty,\qquad k=1,\dots,4,
$$
for any $T^*>0$. Hence alternative \textbf{(B1)} in Corollary~2.5.2 cannot occur, and the local solution extends uniquely to $Q_{T^*}$ for any $T^*>0$. This produces a unique global-in-time weak solution
$\mathbf{u}=(N,T,I,U)$ of \eqref{sys_1}-\eqref{sys_4}. \hfill$\square$

\medskip
The analysis above guarantees that the PDE system \eqref{sys_1}-\eqref{sys_4} is globally well-posed and that the solutions remain nonnegative and bounded on $Q_{T^*}$ for any $T^*>0$. We will next investigate the spatiotemporal behavior of the model and the dosing schedules through numerical simulations. 

\section{CNBE Simulation of the NTIU Model}\label{sec:numerics}

In this section, we investigate the dynamics of \eqref{sys_1}–\eqref{ode_4} through two space-dimensional spatial simulations. We employ a CNBE time discretization in which the linear diffusion–advection operator is advanced implicitly by Crank Nicolson, while the nonlinear reactions are advanced implicitly by Backward Euler. Spatially, we use a conservative finite-volume discretization on a uniform Cartesian grid with centered diffusive fluxes (equivalent to the five-point finite-difference Laplacian in our setting) and first-order upwind advection for the drug, while the nonlinear reactions are advanced implicitly by Backward Euler.

We begin by describing the values of all model parameters, adopting those that are documented in previous tumor–immune models and proposing reasonable estimates for parameters not available in our references. We then describe the computational setting (domain, grid, boundary conditions, initial data). Finally, we conduct three sets of simulations to examine the dynamics of the interacting cell populations under distinct biological scenarios: (i) an immune-free, drug-free NT scheme to establish the baseline dynamics of uncontrolled tumor growth with normal-tissue competition only; (ii) a drug-free NTI scheme to study natural tumor growth with immune response in the absence of chemotherapy; and (iii) a full NTIU system including drug dynamics to evaluate the chemotherapeutic effect on the tumor population.

\subsection{Parameter choices} 
The parameters in our model are primarily adopted from \cite{dePillisRadunskaya2003,DewhirstSecomb2017,FDA,Kuznetsov1994}. In particular, the parameters $(r_1,b_1,s,\rho,\alpha,c_1,c_2,k_1)$ are taken from one of the benchmark tumor-immune models \cite{Kuznetsov1994}. These quantities were originally estimated from BCL$_1$ lymphoma experiments on chimeric mice. The fractional-kill coefficients $a_1, a_2, a_3$ (for immune, tumor, and normal cells, respectively) are calibrated from the work \cite{dePillisRadunskaya2003}, but we increase $a_2$ from $0.3$ to $0.6$ to model a stronger drug effect on tumor cells. The drug decay rate $k_2$ is adapted from U.S. Food and Drug Administration \cite{FDA}, which corresponds to the pharmacokinetic half-life of a chemotherapeutic drug (doxorubicin). In the simulations, we only consider advection in the drug equation ($U$) to reflect directed drug transport through vascular/interstitial pathways \cite{DewhirstSecomb2017}; thus we set $\mathbf{h}_1=\mathbf{h}_2=\mathbf{h}_3=\mathbf{0}$ and take
$\mathbf{h}_4=(10^{-6},10^{-6})$ cm/s.

Several parameters are not directly available from our references and are therefore calibrated. The tumor Allee threshold is set to $A_1 = \frac{0.1}{b_1}$, reflects that a reasonably large tumor population is required to avoid extinction, as small tumor clusters typically do not have enough blood supply, and therefore lead to extinction rather than expansion. The normal tissue Allee threshold is chosen to be $ A_2 = \frac{0.3}{b_2}$. The lower threshold in the tumor population reflects its stronger ability to persist. We choose $r_2 = 0.06\text{day}^{-1}$ and set $b_2 = 1.25\times10^{-9}~\text{cell}^{-1}$, so tumor cells experience faster logistic growth than normal cells. Take $a_0 = 0.1/b_2 = 8 \times10^{7}$ cells to represent the minimal level of normal-tissue density required for drug injection. Set $c_3 = 1\times10^{-11}$ day$^{-1}$\,cell$^{-1}$ and $c_4 = 2\times10^{-11}$ day$^{-1}$\,cell$^{-1}$, the choice $c_4 > c_3$ reflects that tumor cells cause significantly more damage to normal tissue than normal cells do to tumor cells. Following the idea in \cite{DewhirstSecomb2017}, we set
$ D_1 = 1.0\times10^{-6}, D_2 = 8.6\times10^{-5}, D_3 = 1\times10^{-4}, D_4 = 0.086~\text{cm}^2/\text{day}$. These values represent an increasing ability to diffuse from normal to immune cells and a significantly faster diffusion rate for drug particles. All parameters are summarized in the below Table~\ref{tab:params-sim}:

\begin{table}[H]
\centering
\caption{Parameter Values Used in Simulations.}
{\fontsize{8}{10}\selectfont
\setlength{\tabcolsep}{4pt}
\renewcommand{\arraystretch}{1.10}
\begin{tabularx}{\textwidth}{@{} L{0.14\textwidth} L{0.16\textwidth} Y L{0.20\textwidth} @{}}
\toprule
Symbol & Units & Value / Range & Source \\
\midrule
$r_2$ & day$^{-1}$ & $0.06$ & modeling choice \\
$r_1$ & day$^{-1}$ & $0.18$ & \cite{Kuznetsov1994} \\
$b_2$ & cell$^{-1}$ & $1.25\times10^{-9}$ & modeling choice \\
$b_1$ & cell$^{-1}$ & $2.0\times10^{-9}$ & \cite{Kuznetsov1994} \\
$A_2$ & cells & $0.3/b_2$ & modeling choice \\
$A_1$ & cells & $0.1/b_1$ & modeling choice \\
$c_4$ & day$^{-1}$\,cell$^{-1}$ & $2\times10^{-11}$ & modeling choice \\
$c_3$ & day$^{-1}$\,cell$^{-1}$ & $1\times10^{-11}$ & modeling choice \\
$c_2$ & day$^{-1}$\,cell$^{-1}$ & $1.101\times10^{-7}$ & \cite{Kuznetsov1994} \\
$c_1$ & day$^{-1}$\,cell$^{-1}$ & $3.422\times10^{-10}$ & \cite{Kuznetsov1994} \\
$a_1,a_2,a_3$ & -- & $0.2,\,0.6,\,0.1$ & \cite{dePillisRadunskaya2003} \\
$s$ & cells/day & $1.3\times10^{4}$ & \cite{Kuznetsov1994} \\
$\rho$ & day$^{-1}$ & $0.1245$ & \cite{Kuznetsov1994} \\
$\alpha$ & cells & $2.019\times10^{7}$ & \cite{Kuznetsov1994} \\
$k_1$ & day$^{-1}$ & $0.0412$ & \cite{Kuznetsov1994} \\
$k_2$ & day$^{-1}$ & $0.35$ & \cite{FDA} \\
$D_1$ & cm$^2$/day & $1.0\times10^{-6}$ & calibrated based on \cite{DewhirstSecomb2017} \\
$D_2$ & cm$^2$/day & $8.6\times10^{-5}$ & calibrated based on \cite{DewhirstSecomb2017} \\
$D_3$ & cm$^2$/day & $1\times10^{-4}$ & calibrated based on \cite{DewhirstSecomb2017} \\
$D_4$ & cm$^2$/day & $0.086$ & calibrated based on \cite{DewhirstSecomb2017} \\
$a_0$ & cells & $5\times10^{7}$ & modeling choice \\
$\mathbf{h}_4$ & cm/s & $(10^{-6},\,10^{-6})$ & \cite{DewhirstSecomb2017} \\
\bottomrule
\end{tabularx}
}
\label{tab:params-sim}
\end{table}

\subsection{Semi-discrete formulation and CNBE time stepping}\label{subsec:cnbe}
Let $(u_1,u_2,u_3,u_4)=(N,T,I,U)$. Each component satisfies a semi-discrete evolution problem on $t\in[0,T^*]$ of the form
\begin{equation}\label{eq:semi_discrete}
\frac{d}{dt}u_k(t)=\mathcal{L}_k[u_k(t)] + \mathcal{F}_k\big(u_1(t),u_2(t),u_3(t),u_4(t)\big),
\qquad \text{for all } k=1,\dots,4,
\end{equation}
where $\mathcal{L}_k$ is the linear diffusion-advection operator and $\mathcal{F}_k$ collects
the nonlinear reaction functions. In the simulations, advection is nonzero only in the drug equation:
\begin{equation}\label{eq:Lk_choice}
\mathcal{L}_k[u_k]=
\begin{cases}
D_k\,\Delta_h u_k, & k=1,2,3,\\[0.25em]
D_4\,\Delta_h u_4-\mathcal{A}_h u_4, & k=4,
\end{cases}
\end{equation}
where $D_k>0$ are constant diffusion coefficients for $k=1,\cdots,4$, $\Delta_h$ is the five-point discrete Laplacian on a uniform grid, and $\mathcal{A}_h$ is a conservative upwind discretization of advection.

Let $t^n=n\Delta t$ and $u_k^n\approx u_k(t^n)$. The CNBE time discretization for \eqref{eq:semi_discrete}
are:
\begin{equation}\label{eq:cnbe_update}
\frac{u_k^{n+1}-u_k^{n}}{\Delta t}
=\frac12\Big(\mathcal{L}_k[u_k^{n+1}]+\mathcal{L}_k[u_k^{n}]\Big)
+\mathcal{F}_k\big(u_1^{n+1},u_2^{n+1},u_3^{n+1},u_4^{n+1}\big),
\qquad \text{for all } k=1,\dots,4.
\end{equation}

\paragraph{Diffusion Discretization:}
For constant $D_k$ on a uniform grid, we use the standard five-points discretization:
\begin{equation}\label{eq:diff_disc}
(D_k\Delta u_k)_{i,j}\approx
D_k\!(
\frac{(u_k)_{i+1,j}-2(u_k)_{i,j}+(u_k)_{i-1,j}}{\Delta x^2}
+
\frac{(u_k)_{i,j+1}-2(u_k)_{i,j}+(u_k)_{i,j-1}}{\Delta y^2}
).
\end{equation}

\paragraph{Advection Discretization (drug only).}
Let $\mathbf{h}_4=(h_x,h_y)$ with $h_x,h_y\ge 0$. The advection term in conservative form using upwind fluxes gives:
\begin{equation}\label{eq:adv_disc}
(\mathcal{A}_h U)_{i,j}
=\frac{h_x}{\Delta x}\big(U_{i,j}-U_{i-1,j}\big)
+\frac{h_y}{\Delta y}\big(U_{i,j}-U_{i,j-1}\big),
\end{equation}
which is locally conservative.

\subsection{Stability.}
The CNBE scheme treats diffusion (and drug advection) by Crank Nicolson and the nonlinear reactions by Backward Euler. Since Crank Nicolson is A-stable for such linear problems, and Backward Euler is both A-stable and L-stable for the reaction ODE. Therefore the coupled CNBE update is unconditionally stable under these discretizations, see 
\cite{Strikwerda2004}, \cite{HundsdorferVerwer2003}, and \cite{HairerWanner1996}. Although the CNBE scheme is unconditionally stable, the time step is chosen for accuracy. Since the overall temporal accuracy is first order (BE on reactions), we select $\Delta t$ so that (i) advective transport of $U$ does not move more than a fraction of a cell grid per time step, and (ii) there are sufficiently many time steps to capture enough information in each injection duration $\tau$. 

With $\Delta x=\Delta y=0.01$\,cm and $\mathbf{h}_4=(10^{-6},10^{-6})$\,cm/s $=(0.0864,0.0864)$\,cm/day, for advection accuracy, Courant–Friedrichs–Lewy condition requires the Courant number $C \le 1$ (and we take $C=0.7$), where $C$ is calculated by:
$$
C=\frac{|h_x|\Delta t}{\Delta x}+\frac{|h_y|\Delta t}{\Delta y}
=17.28\,\Delta t,
$$ which suggests $\Delta t \le 0.0405$ day. 
To address pulsed dosing accuracy, we require about $10$-$20$ steps per injection duration $\tau$, i.e.
$$
\Delta t \in \Big[\tfrac{\tau}{20},\,\tfrac{\tau}{10}\Big].
$$
Across the tested pulse widths ($\tau\in\{0.30,0.35,0.50\}$ day), this yields $\Delta t\le 0.03$ day. Together, we use $\Delta t=0.025$ day in all simulations.

\subsection{Computational domain, boundary conditions, and initial data}
The computational domain is $\Omega=[0,L_x]\times[0,L_y]$ with $L_x=L_y=1$ cm. 
We use a uniform Cartesian grid with $N_x=N_y=101$ nodes in each direction:
$$
x_i=i\Delta x,\quad y_j=j\Delta y,\qquad 
\text{for } i=0,\dots,N_x-1,\text{ and } j=0,\dots,N_y-1,
$$
where $\Delta x=\frac{L_x}{N_x-1}$ and $\Delta y=\frac{L_y}{N_y-1}$.
Homogeneous Neumann boundary conditions (zero normal flux) are imposed for all components, and are enforced by ghost-point technique. For $k=1,\cdots,4$, we obtain:
$$
(u_k)_{-1,j}=(u_k)_{0,j}, \quad (u_k)_{N_x,j}=(u_k)_{N_x-1,j},\quad
(u_k)_{i,-1}=(u_k)_{i,0},\quad (u_k)_{i,N_y}=(u_k)_{i,N_y-1}
$$

To initialize spatial heterogeneity, we extend the one-dimensional tumor profiles proposed
in \cite{Ansarizadeh2017} to a two-dimensional tissue with three smooth tumor foci locates at $\mathbf{x}_1=(0.35,0.35)$, $\mathbf{x}_2=(0.65,0.35)$, $\mathbf{x}_3=(0.50,0.60)$, and set
$R=0.25$. Define
$$
r_k(x,y)=\sqrt{(x-x_k)^2+(y-y_k)^2},\qquad
G_k(x,y)=
\begin{cases}
\frac{1+\cos\!\big(\pi r_k(x,y)/R\big)}{2}, & r_k(x,y)\le R,\\
0, & r_k(x,y)>R,
\end{cases}
$$
and $B(x,y)=G_1(x,y)+G_2(x,y)+G_3(x,y)$. With carrying capacities
$K_T=1/b_1$ and $K_N=1/b_2$, the initial profile for the tumor and normal density are:
\begin{align}
T_0(x,y)
&=\min\Big\{0.05K_T+0.20K_T\,B(x,y),\ 0.25K_T\Big\},\label{eq:init_T}\\
N_0(x,y)
&=\min\Big\{0.65K_N\big(1-0.10B(x,y)\big),\ 0.75K_N\Big\},\label{eq:init_N}
\end{align}
so that tumor density has peaks where the normal-tissue density has valleys. For immune cells, we mimic the early stage after tumor detection, when immune effector cells begin infiltrating the tumor region from the surrounding tissue. Hence, its density is initially high near the tumor boundary and zero at the tumor center:
$$
I_0(x,y)=\max\Big\{\,3.16\times 10^5\big(0.375-0.235\,\sech^2(r_I(x,y))\big),\ 0\Big\},
\qquad r_I(x,y)=\sqrt{(x-0.50)^2+(y-0.43)^2}.
$$
We assume there is no drug before the treatment, therefore $U_0\equiv 0$.

\subsection{Pulsed dosing protocol and smoothed gating}\label{subsec:dosing}
\begin{figure}[H]
    \centering
    \includegraphics[width=0.5\linewidth]{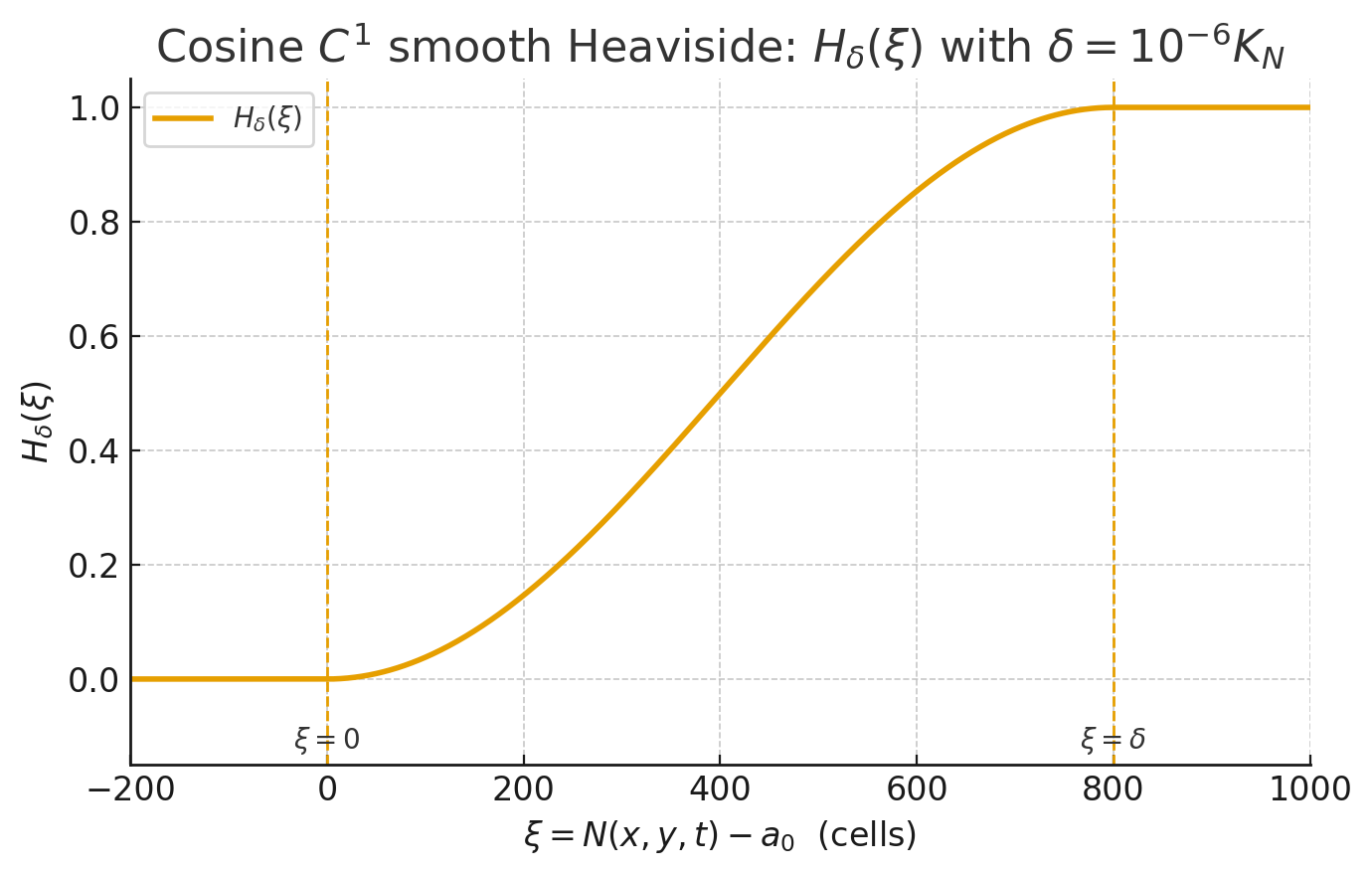}
    \caption{Heaviside Smoothed Function}
    \label{HV_S}
\end{figure}

By the model, treatment drug is being delivered by a pulsed source $v(t)$ that is only active when the normal density stays above a safety threshold $a_0$. The Heaviside gating is replaced by a $C^1$ cosine
smooth-step
\begin{equation}\label{smthed_gate}
\begin{aligned}
H_\delta(\xi)=
\begin{cases}
0, & \xi\le 0,\\
\dfrac{1-\cos(\pi \xi/\delta)}{2}, & 0<\xi<\delta,\\
1, & \xi\ge \delta,
\end{cases}
\end{aligned}
\end{equation}
where $\xi:=N-a_0$, and $\delta=10^{-6}K_N$, shown in the below Figure \ref{HV_S},

\begin{figure}[H]
    \centering
    \includegraphics[width=0.5\linewidth]{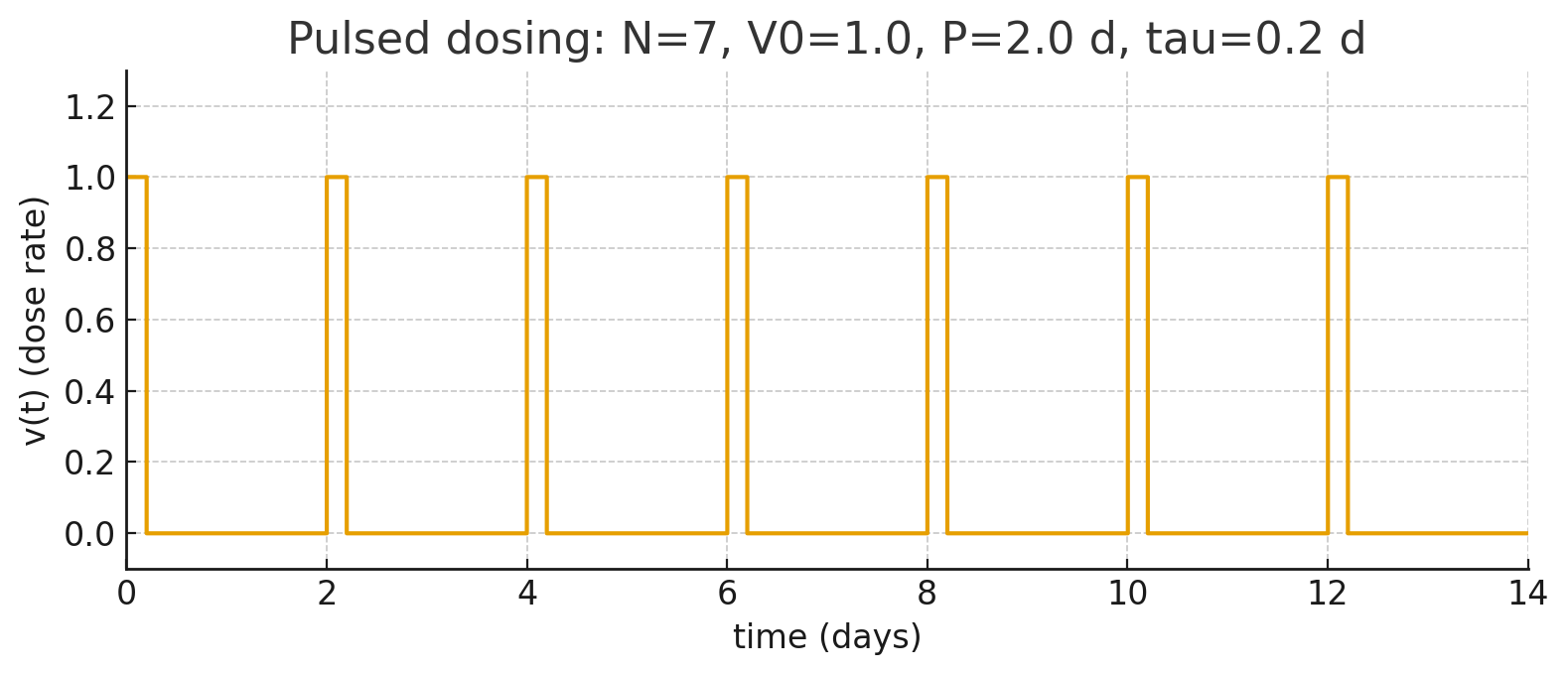}
    \caption{Drug Injection Impulse Function}
    \label{drug_inj}
\end{figure}

Let $\mathbf{1}_{[a,b)}(t)$ denote the indicator on a half-open interval $[a,b)$:
$$
\mathbf{1}_{[a,b)}(t)=
\begin{cases}
1, & a \le t < b,\\
0, & \text{otherwise}.
\end{cases}
$$ Fix $P>0$ be the period between each pulse, $\tau>0$ be the pulse width, let $N^*\in\mathbb{N}$ be the number of pulses, and $V_0>0$ be the per-pulse injection rate. With start
times $t_n=(n-1)P$ for $n= 1,\cdots,N^*$, the dosing source is
$$
v(t)=V_0\sum_{n=1}^{N^*}\mathbf{1}_{[t_n,t_n+\tau)}(t),
$$
so the total delivered dosage over $N^*$ pulses equals $V_0\tau N^*$. For example, Figure \ref{drug_inj} shows the pulsed dosing function $v(t)$ for $N^*=7$ injections: pulses are $P=2$ days apart, each is delivered at rate $V_0=1.0$ lasting $\tau=0.2$ day.

\subsection{Simulations for NT, NTI, and NTIU Subsystems}
We perform two-dimensional simulations for three systems: (i) the NT baseline (normal-tumor competition only), (ii) the NTI subsystem (with immune response included), and (iii) the full NTIU system (with chemotherapy included). In all cases, homogeneous Neumann boundary conditions are imposed and uses the CNBE scheme described in Section 3.2, and parameters follow Table~\ref{tab:params-sim}.

\subsubsection{NT scheme: baseline normal-tumor competition}

The governing equations for the NT model are:
\begin{equation}\label{eq:NT-model}
\begin{aligned}
\partial_t N &= D_1\Delta N
+\underbrace{r_2 N(1-b_2N)\Big(\tfrac{N}{A_2}-1\Big)-c_4TN}_{\mathcal{R}^{NT}_N(N,T)},\\
\partial_t T &= D_2\Delta T
+\underbrace{r_1 T(1-b_1T)\Big(\tfrac{T}{A_1}-1\Big)-c_3TN}_{\mathcal{R}^{NT}_T(N,T)},
\end{aligned}
\end{equation}
subject to homogeneous Neumann boundary conditions:
$$
(D_k\nabla u_k)\!\cdot\!\nu = 0 
\quad\text{on } \partial\Omega, \qquad k=1,2.
$$ The CNBE update equations then read:
\begin{equation}\label{NT_model_D}
\begin{aligned}
\frac{N_{i,j}^{\,n+1}-N_{i,j}^{\,n}}{\Delta t}
&=\frac{D_1}{2}\Big(\Delta_h N_{i,j}^{\,n+1}+\Delta_h N_{i,j}^{\,n}\Big)
+\mathcal{R}_N^{NT}\!\big(N_{i,j}^{\,n+1},T_{i,j}^{\,n+1}\big),\\
\frac{T_{i,j}^{\,n+1}-T_{i,j}^{\,n}}{\Delta t}
&=\frac{D_2}{2}\Big(\Delta_h T_{i,j}^{\,n+1}+\Delta_h T_{i,j}^{\,n}\Big)
+\mathcal{R}_T^{NT}\!\big(N_{i,j}^{\,n+1},T_{i,j}^{\,n+1}\big).
\end{aligned}
\end{equation}

\begin{figure} [H]
    \centering
    \includegraphics[width=1\linewidth]{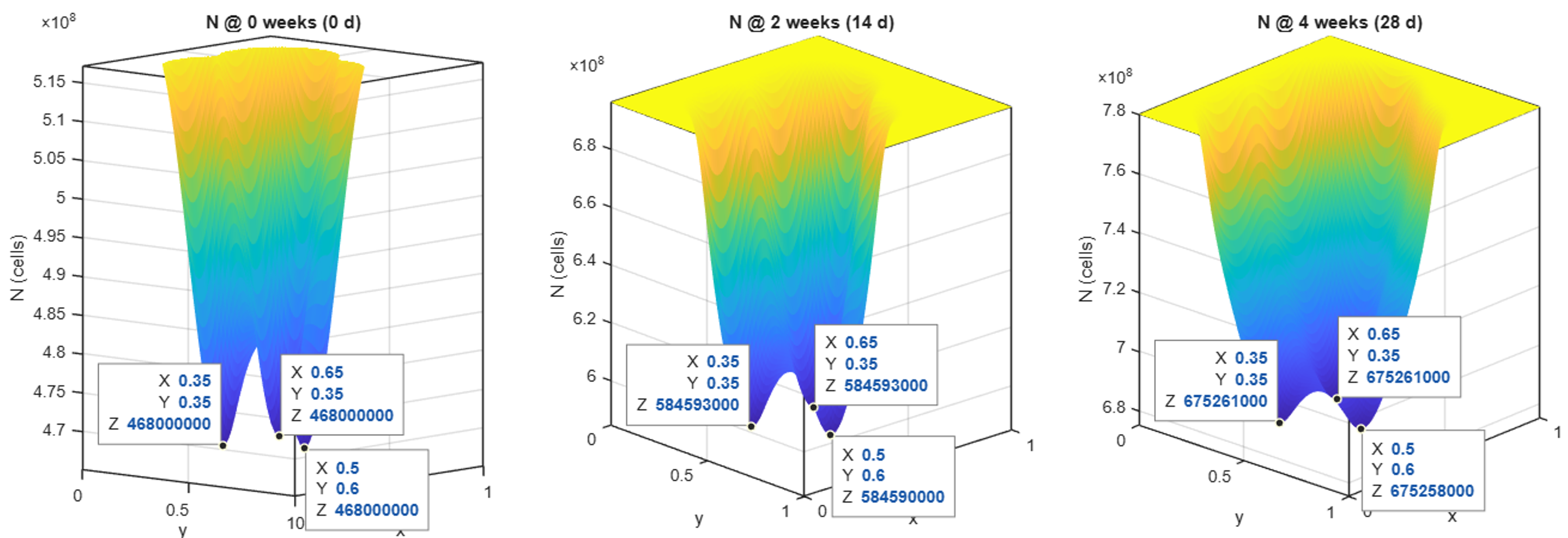}
    \caption{NT Scheme: Normal Cell Density at $t=0$, $2$, $4$ Weeks}
    \label{NT_N}
\end{figure}

\begin{figure}[H]
    \centering
    \includegraphics[width=1\linewidth]{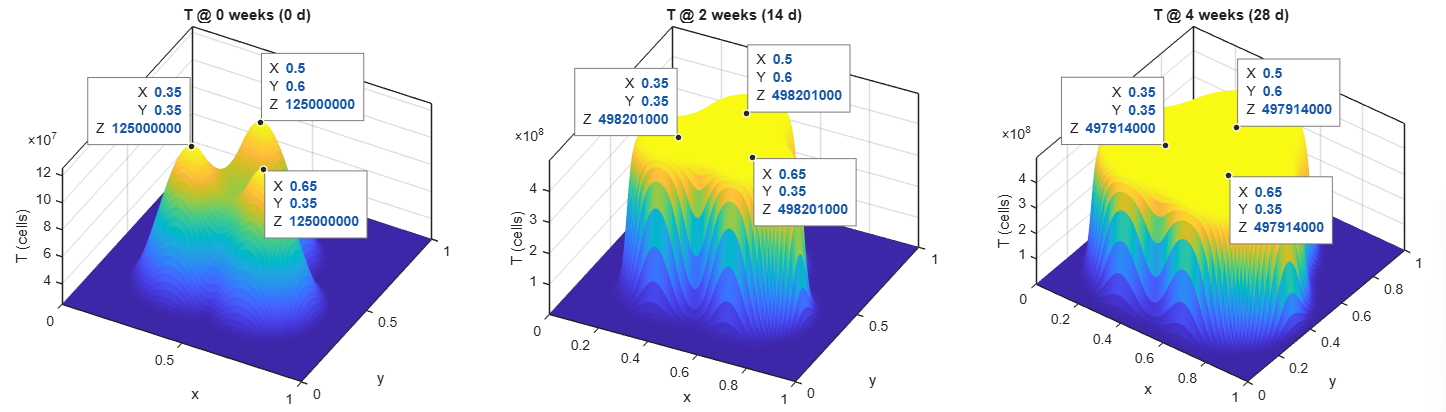}
    \caption{NT Scheme: Tumor Cell Density at $t=0$, $2$, $4$ Weeks}
    \label{NT_T}
\end{figure}

Figures \ref{NT_N} and \ref{NT_T} show the  snapshots of normal and tumor density, respectively at $t=0$, $2$, and $4$ weeks.

Figure \ref{NT_N} shows a mild increase in normal cell density, with the three valleys remaining nearly unchanged in shape. This follows from the small diffusion coefficient $D_1$ (limited mobility), with relatively weak tumor-to-normal competition ($c_4$ small) compared to the logistic with Allee effect term, and the normal density progression is dominated by logistic growth. As a result, the initial pattern is preserved while $N$ slowly drifts toward its carrying capacity $K_N=1/b_2$ over the first four weeks.

Figure \ref{NT_T} shows the temporal evolution of the tumor cell density, it exhibits a steady increase in amplitude and expand spatially: each of the three peaks expands radially due to diffusion, and the three peaks begin to merge into one giant invasive front by the fourth week. The competition term $-c_3TN$ has minimal contribution to tumor density due to small $c_3$. It also shows that without immune and drug effects, the dynamics are dominated by the logistic growth, allowing the tumor density to approach the carrying capacity $K_T=1/b_1$ within the first two weeks. 

\subsubsection{NTI scheme: adding immune response}

The governing equations for the NTI model are:
\begin{equation}\label{eq:NTI-model}
\begin{aligned}
\partial_t N &= D_1\Delta N 
  + \underbrace{r_2 N(1 - b_2 N)(\frac{N}{A_2}-1) - c_4 N T}_{\mathcal{R}_N^{NTI}(N,T,I)},\\
\partial_t T &= D_2\Delta T 
  + \underbrace{r_1 T(1 - b_1 T)(\frac{T}{A_1}-1) - c_3 N T - c_2 I T}_{\mathcal{R}_T^{NTI}(N,T,I)},\\
\partial_t I &= D_3\Delta I 
  + \underbrace{s + \frac{\rho I T}{\alpha + T} - k_1 I - c_1 I T}_{\mathcal{R}_I^{NTI}(N,T,I)},
\end{aligned}
\end{equation}
subject to homogeneous Neumann boundary conditions:
$$
(D_k\nabla u_k)\!\cdot\!\nu = 0 
\quad \text{on } \partial\Omega, \qquad k=1,2,3.
$$ The CNBE update equations then read:
\begin{align}
\frac{N_{i,j}^{\,n+1}-N_{i,j}^{\,n}}{\Delta t}
&=\frac{D_1}{2}\Big(\Delta_h N_{i,j}^{\,n+1}+\Delta_h N_{i,j}^{\,n}\Big)
+\mathcal{R}_N^{NTI}\!\big(N_{i,j}^{\,n+1},T_{i,j}^{\,n+1},I_{i,j}^{\,n+1}\big), \label{eq:NTI-N}\\
\frac{T_{i,j}^{\,n+1}-T_{i,j}^{\,n}}{\Delta t}
&=\frac{D_2}{2}\Big(\Delta_h T_{i,j}^{\,n+1}+\Delta_h T_{i,j}^{\,n}\Big)
+\mathcal{R}_T^{NTI}\!\big(N_{i,j}^{\,n+1},T_{i,j}^{\,n+1},I_{i,j}^{\,n+1}\big), \label{eq:NTI-T}\\
\frac{I_{i,j}^{\,n+1}-I_{i,j}^{\,n}}{\Delta t}
&=\frac{D_3}{2}\Big(\Delta_h I_{i,j}^{\,n+1}+\Delta_h I_{i,j}^{\,n}\Big)
+\mathcal{R}_I^{NTI}\!\big(N_{i,j}^{\,n+1},T_{i,j}^{\,n+1},I_{i,j}^{\,n+1}\big), \label{eq:NTI-I}
\end{align}

\begin{figure} [H]
    \centering
    \includegraphics[width=1\linewidth]{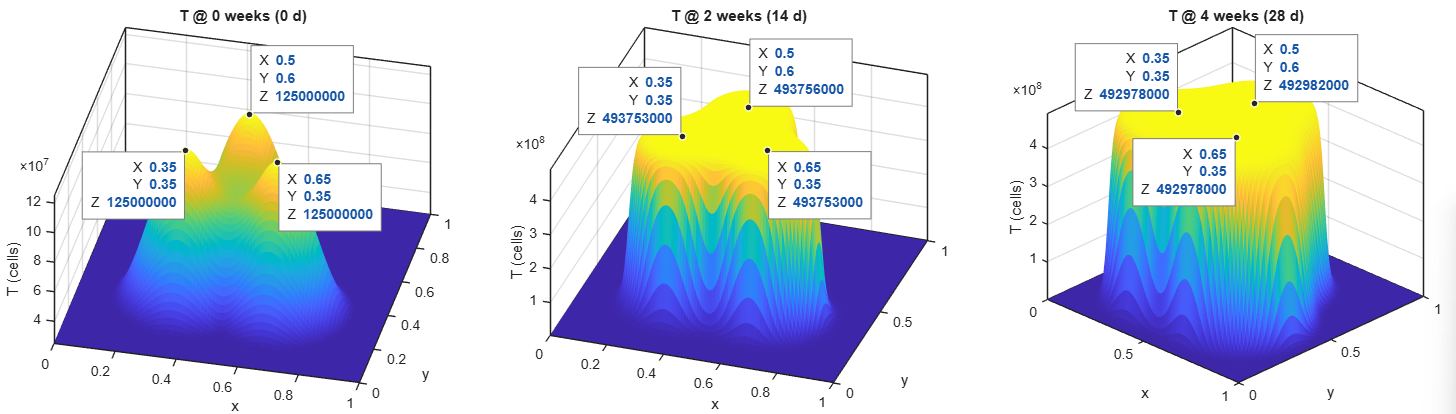}
    \caption{NTI Scheme: Tumor Cell Density at $t=0$, $2$, $4$ Weeks}
    \label{NTI_T}
\end{figure}

\begin{figure} [H]
    \centering
    \includegraphics[width=1\linewidth]{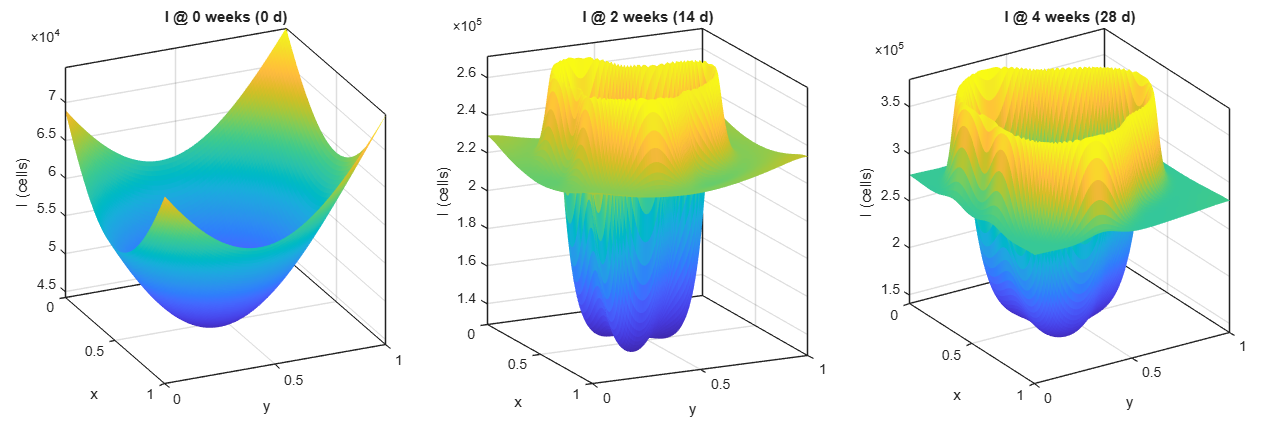}
    \caption{NTI Scheme: Immune Cell Density at $t=0$, $2$, $4$ Weeks}
    \label{NTI_I}
\end{figure}

Figures \ref{NTI_N}, \ref{NTI_T}, and \ref{NTI_I} show the results for the NTI scheme:

\begin{figure} [H]
    \centering
    \includegraphics[width=1\linewidth]{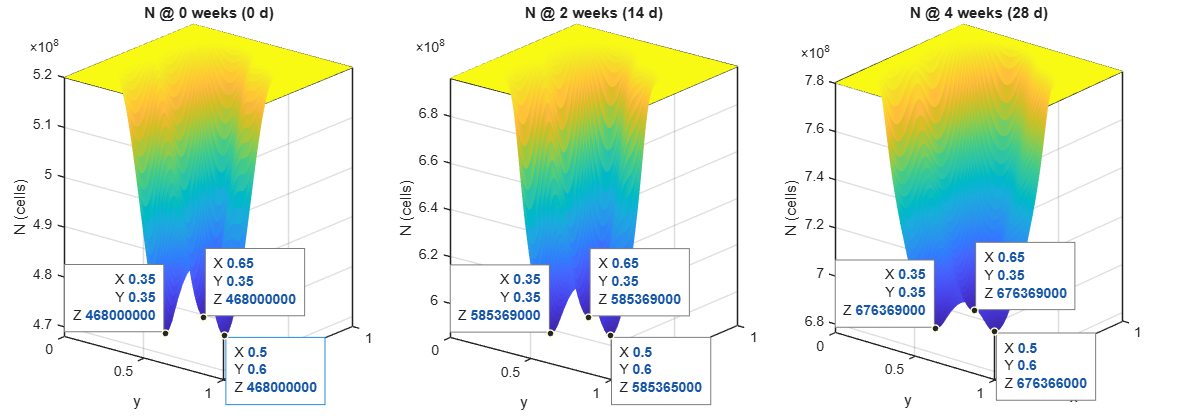}
    \caption{NTI Scheme: Normal Cell Density at $t=0$, $2$, $4$ Weeks}
    \label{NTI_N}
\end{figure}

The most notable feature of the NTI scheme is the evolution of the immune density. Figure \ref{NTI_I} shows that $I$ increases in amplitude over $0-4$ weeks, indicating that the net local gain from the source and proliferation term $s+\rho IT/(\alpha+T)$ outweighs the losses $k_1I+c_1IT$. It is also observed that by week 2, $I$ organizes into an annular band surrounding the tumor peaks, and the ring thickens and its inner edge advances inward. From week 2 to week 4, the band continues to penetrate toward the tumor interior but shifts slightly backward the tumor center. This is because near the inner edge, where $T$ becomes large, the immune reaction therefore becomes negative and $I$ decays, but along the outer edge, $T$ is smaller, and the immune population experience positive net growth. Overall, observation indicates a active immune response near tumor regions and, due to the larger mobility $D_3$, immune effector cells penetrate toward the tumor interior.

Figure \ref{NTI_T} shows that, despite immune effects, the tumor still grows rapidly toward its carrying capacity $K_T$ and forms a single broad invasive front by week 4. Relative to the NT case, the peak values at the initial tumor foci are reduced, but only by a small margin, and the plateau remains close to $K_T$. Which indicates that in our parameter setting, the immune response only tempers but does not control tumor expansion.

Comparing Figures \ref{NTI_N} and \ref{NT_N}, normal tissue shows no visible change in shape and no substantial change in the depths of the three valleys. This is expected because the $N$-equation is identical in the NT and NTI subsystems, so the logistic (with Allee effect) dynamics for $N$ are unchanged. Moreover, the tumor trajectory in the NTI simulation remains close to that in the NT case, so the coupling term $-c_4NT$ has a similar magnitude in both subsystems.

\subsubsection{NTIU scheme: full system with pulsed chemotherapy}\label{subsubsec:NTIU}

Finally, We simulate the full NTIU system with advection included only in the drug equation.
The governing PDEs are
\begin{equation}\label{eq:NTIU-model}
\begin{aligned}
\partial_t N &= D_1\Delta N
+\mathcal{R}^{NTIU}_N(N,T,I,U),\\
\partial_t T &= D_2\Delta T
+\mathcal{R}^{NTIU}_T(N,T,I,U),\\
\partial_t I &= D_3\Delta I
+\mathcal{R}^{NTIU}_I(N,T,I,U),\\
\partial_t U &= D_4\Delta U - \nabla\!\cdot(\mathbf{h}_4U)
+\mathcal{R}^{NTIU}_U(N,T,I,U),
\end{aligned}
\end{equation}
where the reaction terms $\mathcal{R}^{NTIU}_N, \mathcal{R}^{NTIU}_T, \mathcal{R}^{NTIU}_I$ and  $\mathcal{R}^{NTIU}_U$ follow those in \eqref{ode_1}-\eqref{ode_4} with the smoothed gate \eqref{smthed_gate} replacing $H_\delta(N-a_0)$ and the pulsed input $v(t)$. subject to homogeneous Neumann boundary conditions:
$$
(D_k\nabla u_k)\cdot\nu = 0 \ (k=1,2,3),\qquad
(D_4\nabla U - \mathbf{h}_4U)\cdot\nu = 0 \quad \text{on }\partial\Omega.
$$

The CNBE update equations then read:
\begin{align}
\frac{N_{i,j}^{n+1}-N_{i,j}^{n}}{\Delta t}
&=\frac{D_1}{2}\big(\Delta_hN_{i,j}^{n+1}+\Delta_hN_{i,j}^{n}\big)
+\mathcal{R}^{NTIU}_N\big(N_{i,j}^{\,n+1},T_{i,j}^{\,n+1},I_{i,j}^{\,n+1},U_{i,j}^{\,n+1}\big),\label{eq:NTIU-cnbe-N}\\
\frac{T_{i,j}^{n+1}-T_{i,j}^{n}}{\Delta t}
&=\frac{D_2}{2}\big(\Delta_hT_{i,j}^{n+1}+\Delta_hT_{i,j}^{n}\big)
+\mathcal{R}^{NTIU}_T\big(N_{i,j}^{\,n+1},T_{i,j}^{\,n+1},I_{i,j}^{\,n+1},U_{i,j}^{\,n+1}\big),\label{eq:NTIU-cnbe-T}\\
\frac{I_{i,j}^{n+1}-I_{i,j}^{n}}{\Delta t}
&=\frac{D_3}{2}\big(\Delta_hI_{i,j}^{n+1}+\Delta_hI_{i,j}^{n}\big)
+\mathcal{R}^{NTIU}_I\big(N_{i,j}^{\,n+1},T_{i,j}^{\,n+1},I_{i,j}^{\,n+1},U_{i,j}^{\,n+1}\big),\label{eq:NTIU-cnbe-I}\\
\frac{U_{i,j}^{n+1}-U_{i,j}^{n}}{\Delta t}
&=\frac{D_4}{2}\big(\Delta_hU_{i,j}^{n+1}+\Delta_hU_{i,j}^{n}\big)
-\frac12\big(\mathcal{A}_hU_{i,j}^{n+1}+\mathcal{A}_hU_{i,j}^{n}\big)\\
&+\mathcal{R}^{NTIU}_U\big(N_{i,j}^{\,n+1},T_{i,j}^{\,n+1},I_{i,j}^{\,n+1},U_{i,j}^{\,n+1}\big).\label{eq:NTIU-cnbe-U}
\end{align}

For the full NTIU system, we compare four pulsed-therapy cases listed in Table \ref{tab:NTIU-cases}. Cases 1-2 share the same total dosage, while Cases 3-4 share a larger (but same) total dosage.

\begin{table}[H]
\centering
\caption{NTIU dosing cases.}
{\fontsize{8}{10}\selectfont
\setlength{\tabcolsep}{10pt}
\renewcommand{\arraystretch}{1.15}
\begin{tabular}{c c c c}
\toprule
Case & $V_0$ & $\tau$ (day) & $N^*$ \\
\midrule
1 & 1.00 & 0.30 & 7 \\
2 & 0.60 & 0.50 & 7 \\
3 & 3.00 & 0.35 & 7 \\
4 & 2.10 & 0.35 & 10 \\
\bottomrule
\end{tabular}
}
\label{tab:NTIU-cases}
\end{table}

\begin{figure} [H]
    \centering
    \includegraphics[width=1\linewidth]{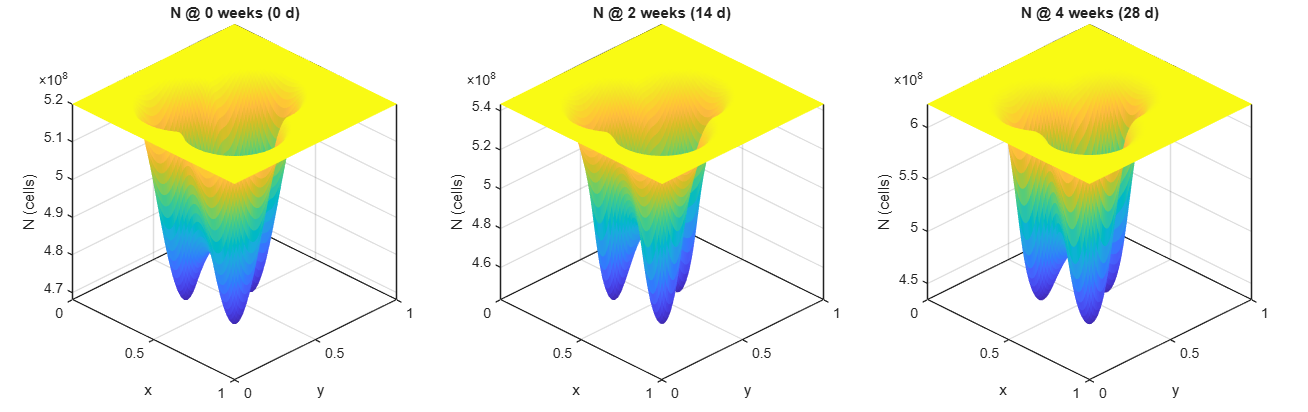}
    \caption{NTIU Scheme(Case $1$): Normal Cell Density at $t=0$, $2$, $4$ Weeks}
    \label{CASE1_N}
\end{figure}

\begin{figure}[H]
    \centering
    \includegraphics[width=1\linewidth]{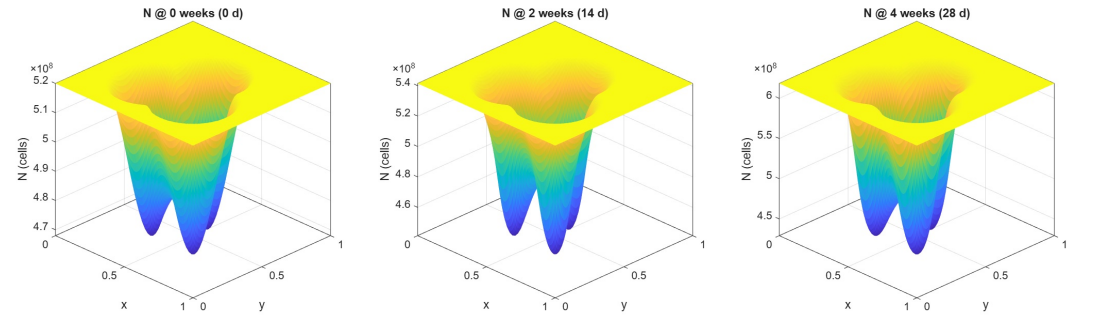}
    \caption{NTIU Scheme(Case $2$): Normal Cell Density at $t=0$, $2$, $4$ Weeks}
    \label{CASE2_N}
\end{figure}

\begin{figure}[H]
    \centering
    \includegraphics[width=1\linewidth]{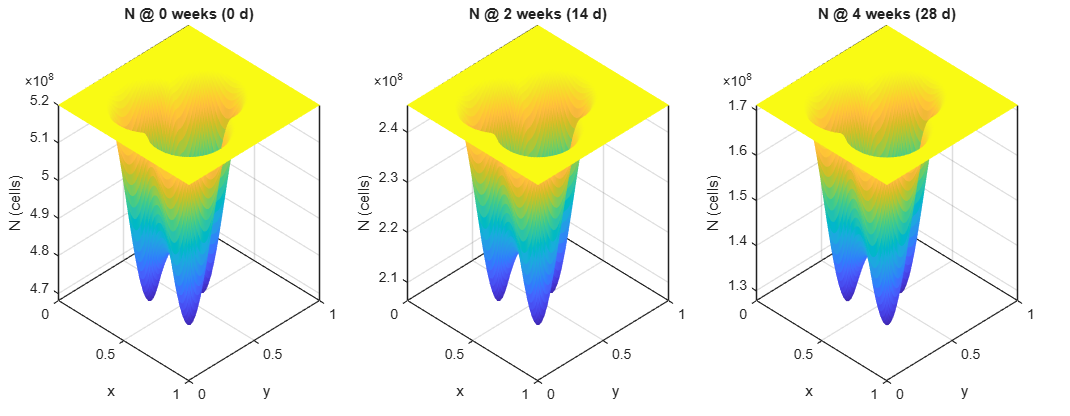}
    \caption{NTIU Scheme(Case $3$): Normal Cell Density at $t=0$, $2$, $4$ Weeks}
    \label{CASE3_N}
\end{figure}

\begin{figure}[H]
    \centering
    \includegraphics[width=1\linewidth]{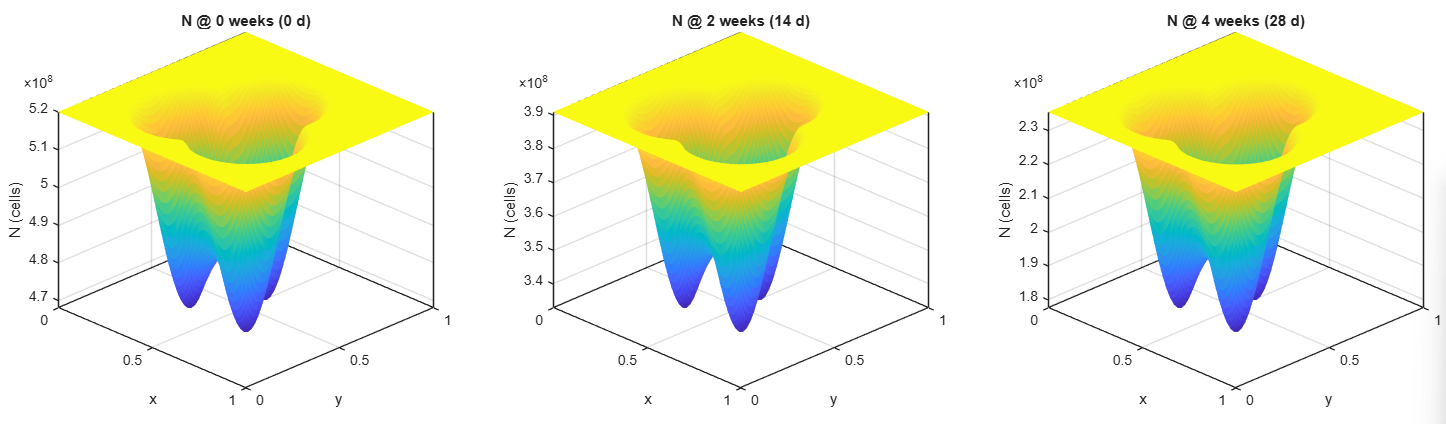}
    \caption{NTIU Scheme(Case $4$): Normal Cell Density at $t=0$, $2$, $4$ Weeks}
    \label{CASE4_N}
\end{figure}

\paragraph{Results for Normal Cell Density:}
Figure \ref{CASE1_N}–\ref{CASE4_N} show the normal-cell density for Cases 1-4 at $t=0$, $2$, and $4$ weeks.

Across all cases, the normal cell density retains its initial spatial profile because diffusion is weak (small $D_1$). Cases 1 and 2 have similar amplitudes, whereas cases 3 and 4 are significantly lower, with case 3 being the lowest.

\begin{figure}[H]
    \centering
    \includegraphics[width=1\linewidth]{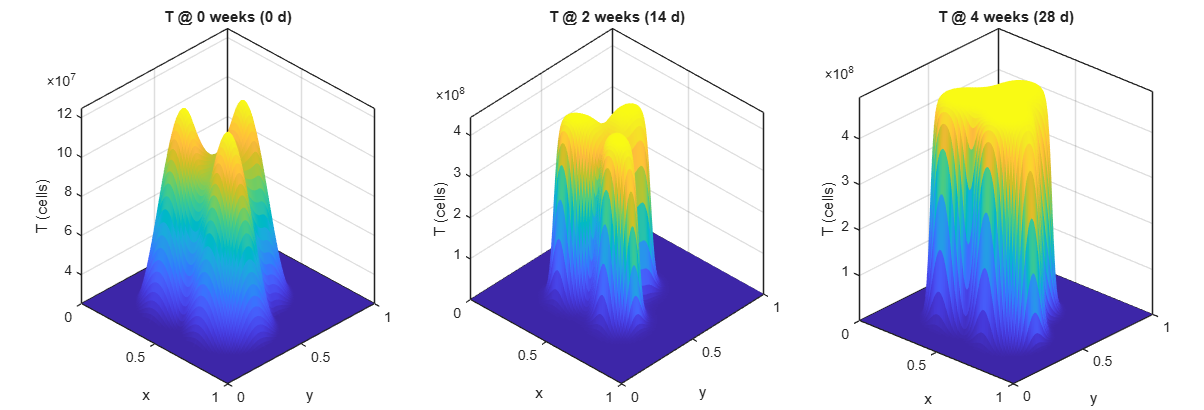}
    \caption{NTIU Scheme(Case $1$): Tumor Cell Density at $t=0$, $2$, $4$ Weeks}
    \label{CASE1_T}
\end{figure}

\begin{figure}[H]
    \centering
    \includegraphics[width=1\linewidth]{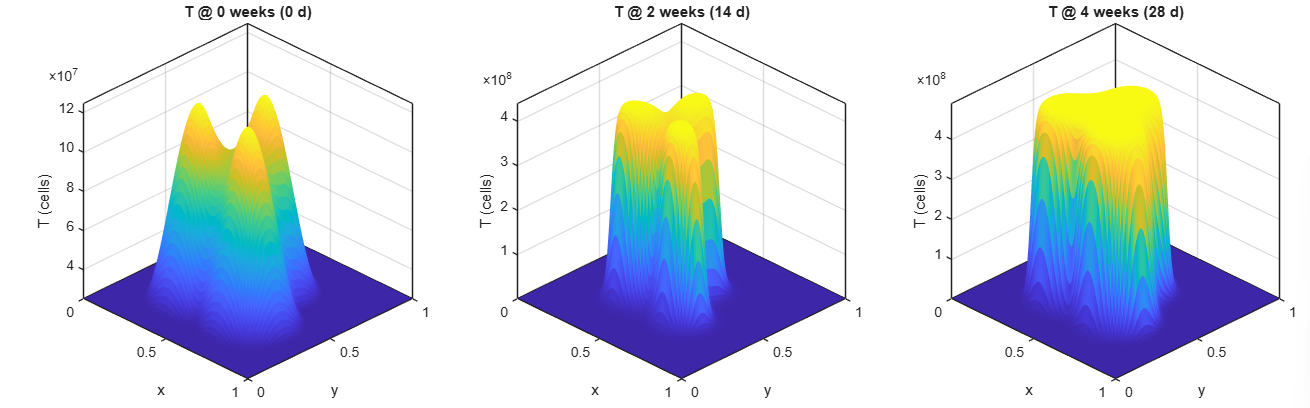}
    \caption{NTIU Scheme(Case $2$): Tumor Cell Density at $t=0$, $2$, $4$ Weeks}
    \label{CASE2_T}
\end{figure}

\begin{figure}[H]
    \centering
    \includegraphics[width=1\linewidth]{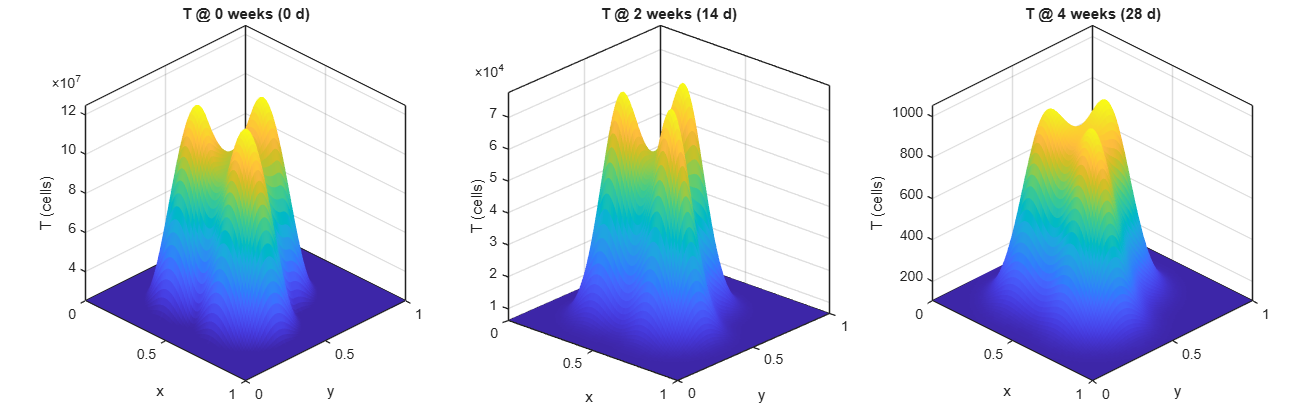}
    \caption{NTIU Scheme(Case $3$): Tumor Cell Density at $t=0$, $2$, $4$ Weeks}
    \label{CASE3_T}
\end{figure}

\begin{figure}[H]
    \centering
    \includegraphics[width=1\linewidth]{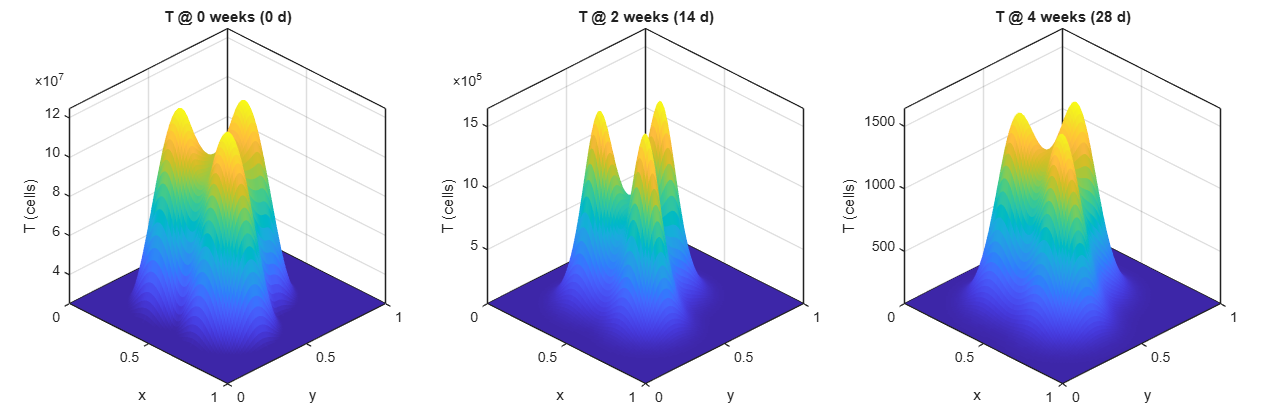}
    \caption{NTIU Scheme(Case $4$): Tumor Cell Density at $t=0$, $2$, $4$ Weeks}
    \label{CASE4_T}
\end{figure}

\paragraph{Results for Tumor Cell Density:}
Figure \ref{CASE1_T}-\ref{CASE4_T} show the tumor-cell density for Cases 1-4 at $t=0$, $2$, and $4$ weeks.

In Cases 1-2, the tumor-density surfaces are similar in both shape and amplitude at weeks $2$ and $4$. Relative to the NT and NTI subsystems, both cases show a clear retreat of the invasive front at both time-points, while the peak amplitudes are only mildly reduced. Cases 3-4 yield an even stronger overall tumor suppression than Cases 1-2. In particular, Case 3 produces the largest reduction in peak amplitude at both weeks $2$ and $4$. Case 4 achieves the strongest contraction of the invasive front at both weeks $2$ and $4$; although its peak remains higher than Case 3 at week $2$, by week $4$ the two cases show comparable peak amplitudes.

\begin{figure}[H]
    \centering
    \includegraphics[width=1\linewidth]{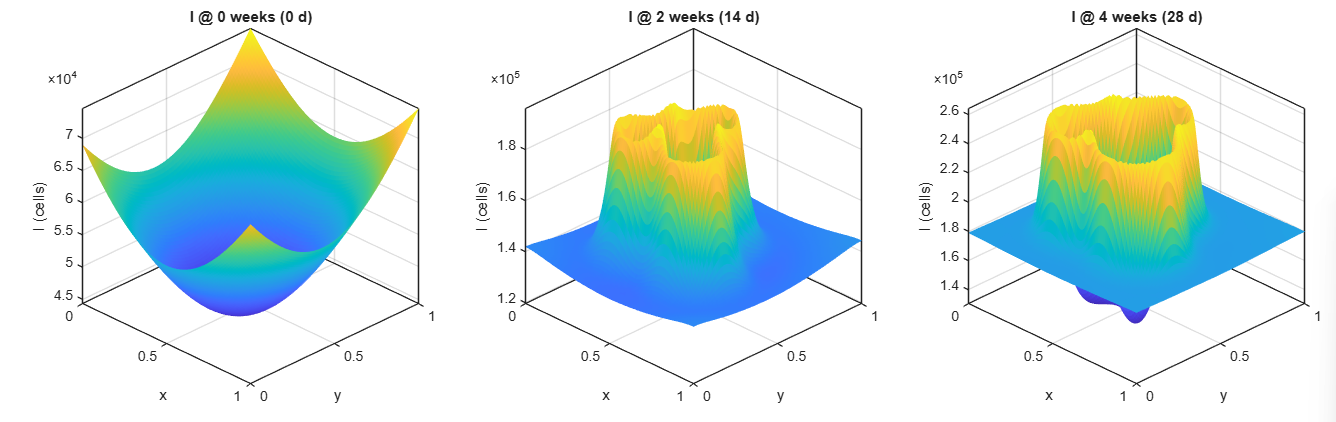}
    \caption{NTIU Scheme(Case $1$): Immune Cell Density at $t=0$, $2$, $4$ Weeks}
    \label{CASE1_I}
\end{figure}

\begin{figure}[H]
    \centering
    \includegraphics[width=1\linewidth]{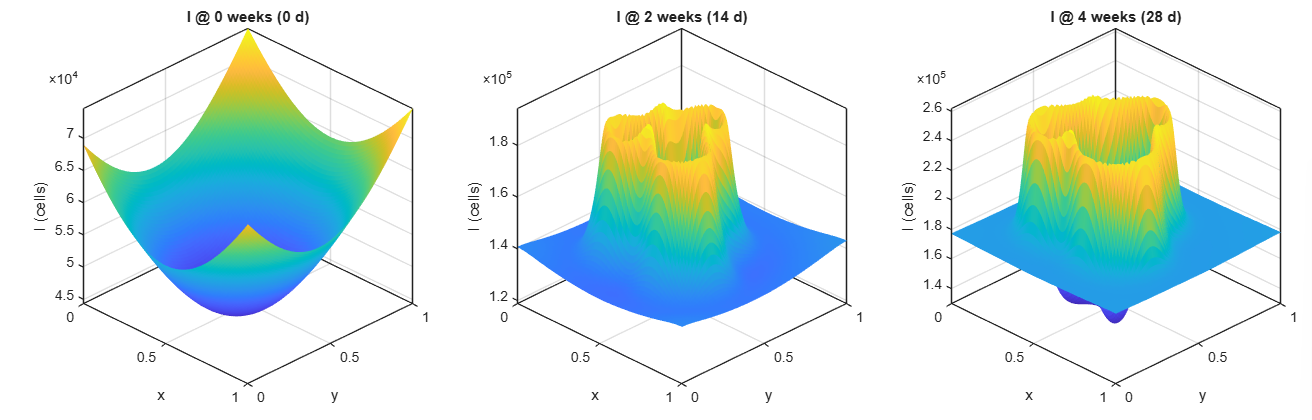}
    \caption{NTIU Scheme(Case $2$): Immune Cell Density at $t=0$, $2$, $4$ Weeks}
    \label{CASE2_I}
\end{figure}

\begin{figure}[H]
    \centering
    \includegraphics[width=1\linewidth]{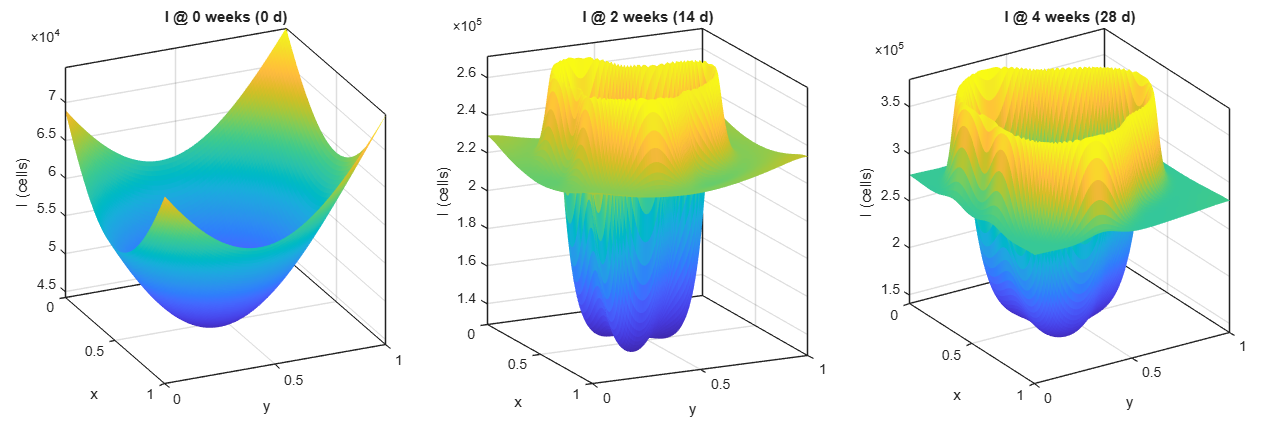}
    \caption{NTIU Scheme(Case $3$): Immune Cell Density at $t=0$, $2$, $4$ Weeks}
    \label{CASE3_I}
\end{figure}

\begin{figure}[H]
    \centering
    \includegraphics[width=1\linewidth]{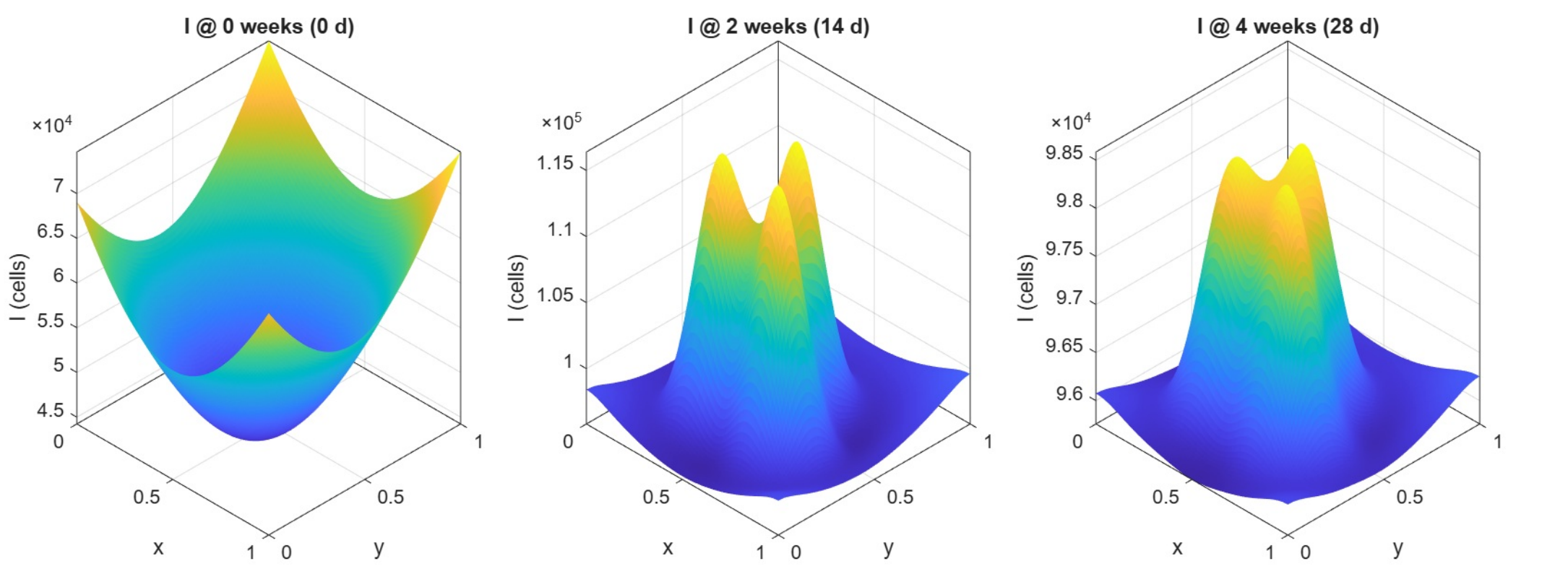}
    \caption{NTIU Scheme(Case $4$): Immune Cell Density at $t=0$, $2$, $4$ Weeks}
    \label{CASE4_I}
\end{figure}

\paragraph{Results for Immune Cell Density:}
Figure  \ref{CASE1_I}-\ref{CASE4_I} show the immune-cell density for Cases 1-4 at $t=0$, $2$, and $4$ weeks.

From the simulations, the immune density surfaces in cases 1-3 have a similar shape: an annular band encircling the tumor peaks. Cases $1$ and $2$ have comparable amplitudes, whereas case $3$ shows a noticeably higher crest at the top of the band. In case $4$, the immune density forms three peaks matched with the tumor foci at weeks $2$ and $4$; however, these peak amplitudes are lower than the annular-band crest observed in case $3$.

\begin{figure}[H]
    \centering
    \includegraphics[width=1\linewidth]{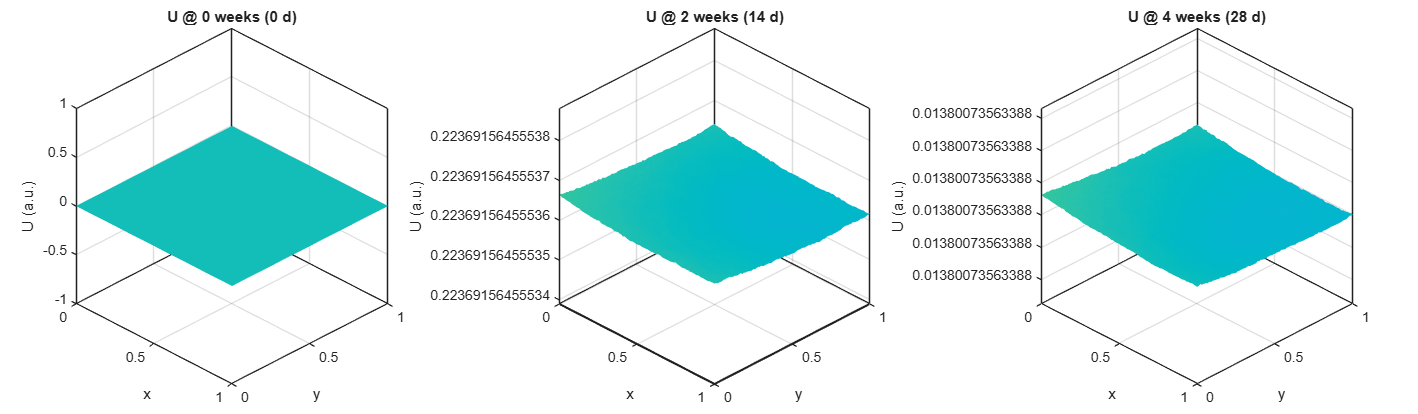}
    \caption{NTIU Scheme(Case $1$): Drug Density at $t=0$, $2$, $4$ Weeks}
    \label{CASE1_U}
\end{figure}

\begin{figure}[H]
    \centering
    \includegraphics[width=1\linewidth]{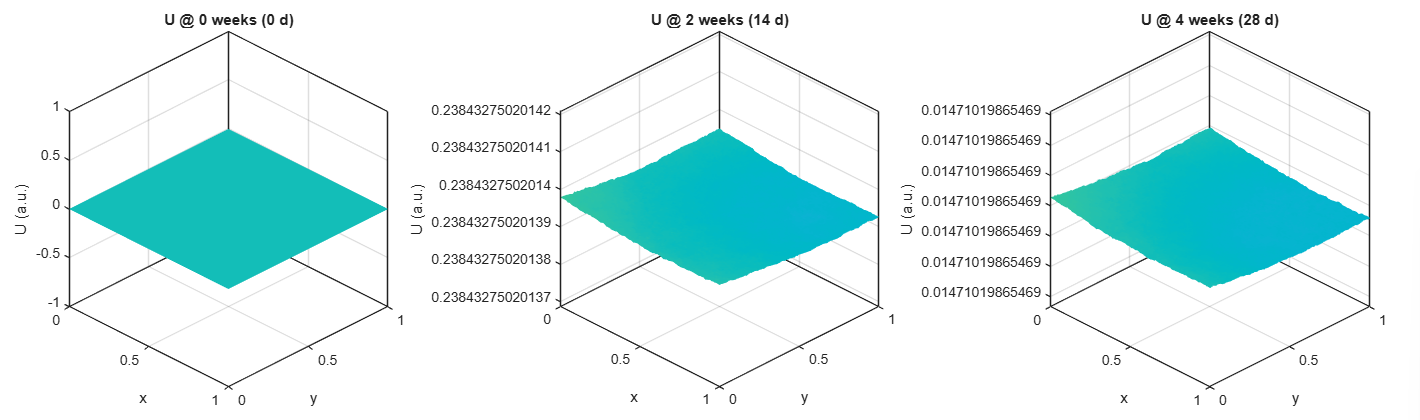}
    \caption{NTIU Scheme(Case $2$): Drug Density at $t=0$, $2$, $4$ Weeks}
    \label{CASE2_U}
\end{figure}

\begin{figure}[H]
    \centering
    \includegraphics[width=1\linewidth]{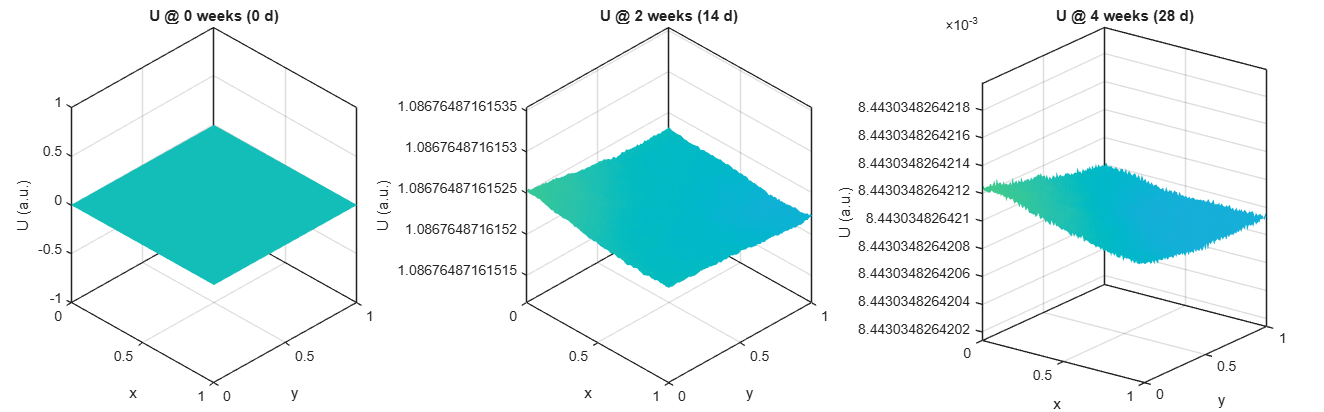}
    \caption{NTIU Scheme(Case $3$): Drug Density at $t=0$, $2$, $4$ Weeks}
    \label{CASE3_U}
\end{figure}

\begin{figure}[H]
    \centering
    \includegraphics[width=1\linewidth]{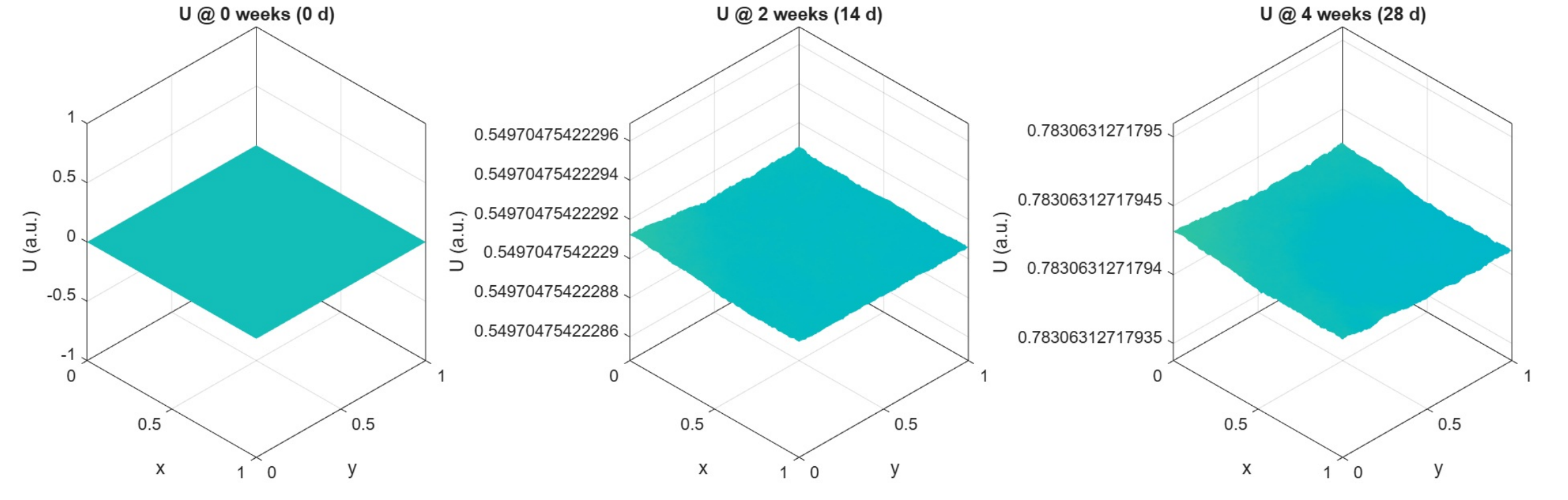}
    \caption{NTIU Scheme(Case $4$): Drug Density at $t=0$, $2$, $4$ Weeks}
    \label{CASE4_U}
\end{figure}

\paragraph{Results for Drug Density:}
Figure \ref{CASE1_U}-\ref{CASE4_U} show the result for drug concentration from case 1-4, at $0$, $2$, and $4$ weeks.

Across all cases, the drug concentration surfaces are nearly planar, indicating a broadly uniform spatial distribution.

\paragraph{Results and Discussion}
In the NT scheme we examine the normal–tumor dynamics without immune or drug. The tumor density rises rapidly toward its carrying capacity and, by week $2$, the three peaks merge into a single broad invasive front, and by week $4$ this front has expanded further. With immune response, NTI scheme shows the immune density forms an annular band around the tumor foci. Its amplitude increases from week $2$ to week $4$ and the band thickens. It is also observed that immune proliferation reduces the tumor only slightly and does not overturn the dominated tumor growth. The normal tissue is affected minimally.

For the NTIU scheme, we test four treatment cases. In cases 1-2 (lower–dosage), the total drug amount is smaller than in cases 3-4. Drug delivery markedly suppresses the tumor front by week $2$ in all cases, in particular, the front remains smaller than in NT or NTI scheme across all cases at both week $2$ and $4$. In cases 1-2, the peak tumor amplitude is only mildly reduced, and both normal and immune densities experience small decreases. It is observed that when the total drug and the number of injections are held fixed, varying the per–pulse rate versus duration yields similar overall outcomes.

By contrast, cases 3-4 use a higher overall dosage, while holding the total dosage equal between the two cases, case 3 delivers it in $7$ higher–rate pulses, whereas case 4 uses $10$ lower–rate pulses. Case 3 produces the largest early reduction of the tumor peak amplitude at week $2$, while case 4 achieves a stronger early contraction of the invasive front at both week $2$ and $4$. By week $4$ the overall tumor reduction is comparable in the two setting. In the higher-dosage schedules (cases 3-4), the larger cumulative exposure to $U$ also suppresses the immune population more strongly, reflecting the immunotoxic side effects of chemotherapy encoded by the $a_1\big(1-e^{-U}\big)I$ term. Notably, case $4$ causes less collateral damage in normal density, suggesting that distributing the same dosage in smaller but more pulses can better spare normal tissue while still controlling spread of tumor to the similar level as in case $3$ at the end of the treartment cycle, i.e., a longer treatment cycle (more, gentler injections) may be preferable when normal–tissue preservation is considered.

Together, these results suggest that for a fixed total amount of drug, distributing the dosage into more frequent, lower-intensity pulses can better preserve normal tissue while achieving tumor control comparable to that of fewer, more aggressive pulses.

\section{Conclusion and Future Work}
In this paper, we proposed and analyzed a four-species reaction-diffusion-advection semilinear parabolic model for normal tissue $N$, tumor cells $T$, immune effector cells $I$, and chemotherapeutic drug $U$ under homogeneous Neumann boundary conditions. Two features distinguish the work: (i) strong Allee effects are incorporated in both  $N$ and $T$ growth terms to encode extinction below critical thresholds, and (ii) each equation is posed with a general second-order elliptic operator in divergence form, providing a flexible framework for heterogeneous tissue transport and directed drug delivery. On the analytical side, under minimal structural assumptions we established global well-posedness of nonnegative weak solutions and derived uniform $L^\infty(Q_{T^*})$ bounds on any finite time interval. In particular, boundedness on $[0,T^*]$ rules out finite-time blow-up and yields a unique global solution via a standard blow-up alternative. On the computational side, we developed a conservative CNBE scheme, treating diffusion-advection implicitly and nonlinear reactions by Backward Euler, together with ghost-point Neumann boundaries. The simulations illustrate a consistent progression from baseline tumor invasion (NT), to the immune effect subsystem without chemotherapy (NTI), and then to tumor suppression under pulsed dosing (NTIU). For a fixed total dosage per cycle, we observe a scheduling trade-off: fewer, higher-rate pulses produce stronger early peak knock-down, while more frequent, gentler pulses provide comparable control of the invasive front with improved normal-tissue preservation. We would like to mention a recent paper \cite{MTY2024} that deals with a similar model with injection of a drug through the tissue boundary.

Several extensions are natural. Analytically, it would be of interest to upgrade weak solutions to stronger regularity classes under additional smoothness assumptions and to study long-time dynamics. From a modeling and numerical perspective, the sharp strong-Allee factor could be replaced by a smoothed Allee function to reduce threshold-induced numerical difficulties.

\section{Appendix A}

\subsection{Proof of Theorem 2.5.1 }

\noindent\textit{Proof.}
Given a sequence $\{\mathbf{w}^j\}_{j=1}^\infty$ converging to $\mathbf{w}\in K$ in
$$
\|\mathbf{w}^j-\mathbf{w}\|_{L^{\infty}(Q_{T^*};\mathbb{R}^4)} \rightarrow 0 \quad \text{as } j\rightarrow \infty .
$$
Set $\mathbf{u}^{*j}=M[\mathbf{w}^j]$, $\mathbf{u}^*=M[\mathbf{w}]$, and $\mathbf{v}^j=\mathbf{u}^{*j}-\mathbf{u}^*$.
For each component $k=1,\dots,4$, $v_k^j$ solves
\begin{align}
&\partial_t v^j_k+L_k[v^j_k] \,=\, F_k(x,t,\mathbf{w}^j)\,-\,F_k(x,t,\mathbf{w}), && (x, t) \in Q_{T^*},  \label{con_2}\\
&(D_k(x,t)\,\nabla v^j_k)\cdot\nu\,=\,0, && (x, t) \in S_{T^*}, \label{con_3}\\
&v^j_k(x, 0) \,=\, 0, && x \in \Omega. \label{con_4}
\end{align}

First, $F_k(x,t,\mathbf{w})$ is locally Lipschitz in $\mathbf{w}$ on $K$ for each $k=1,\dots,4$; let $L$ be a Lipschitz constant such that
\begin{equation}\label{lip_L}
|F_k(x,t,\mathbf{w}^j)-F_k(x,t,\mathbf{w})|\le L\,|\mathbf{w}^j-\mathbf{w}|.
\end{equation}
Multiply \eqref{con_2} by $2v_k^j$, integrate over $\Omega$, sum over $k=1,\dots,4$, and integrate by parts using \eqref{con_3} and \textbf{H(1)}.
Using \eqref{lip_L} and Young's inequality (with $\epsilon_1,\epsilon_2>0$), we obtain
\begin{align}\label{eq:cont3}
\frac{d}{dt}\sum_{k=1}^4\int_{\Omega}\big (v^j_k(\cdot,t)\big)^2\,dx
\le C_1\sum_{k=1}^4\int_{\Omega}\big (v^j_k(\cdot,t)\big)^2\,dx
+ C_2\sum_{\ell=1}^4\int_{\Omega}\big(w_\ell^j-w_\ell\big)^2\,dx,
\end{align}
where $C_1=C_1(d_0,h^*,\epsilon_1,\epsilon_2)$ and $C_2=C_2(L,\epsilon_2)$.
Choose $\epsilon_1=\frac{d_0}{2h^*}$ and $\epsilon_2=1$, so $C_1=\frac{(h^*)^2}{d_0}+1>0$ and $C_2=4L^2$.
Applying Grönwall's inequality to \eqref{eq:cont3} and using $v_k^j(\cdot,0)=0$ yields
\begin{equation}\label{eq:cont4}
\sup_{0\le t\le T^*}\sum_{k=1}^4\|v_k^j(\cdot,t)\|_{L^2(\Omega)}^2
\le C(T^*)\,\|\mathbf{w}^j-\mathbf{w}\|_{L^2(Q_{T^*})}^2 \;\rightarrow\;0 .
\end{equation}

By the definition of $M$ and the H\"older estimate \eqref{DN_est}, $\{\mathbf{u}^{*j}\}$ is uniformly bounded in $C^{\beta,\beta/2}(\overline{Q_{T^*}};\mathbb{R}^4)$, the family $\{\mathbf{v}^j\}$ is uniformly bounded and equicontinuous on $\overline{Q_{T^*}}$. Thus, by Arzel\`a--Ascoli, any subsequence of $\{\mathbf{v}^j\}$ has a further subsequence converging in $C(\overline{Q_{T^*}};\mathbb{R}^4)$ to some limit $\mathbf{v}$. The $L^2$ convergence in \eqref{eq:cont4} forces $\mathbf{v}\equiv \mathbf{0}$, hence the whole sequence satisfies
$$
\|\mathbf{v}^j\|_{L^\infty(Q_{T^*};\mathbb{R}^4)}\rightarrow 0.
$$
Therefore $\mathbf{u}^{*j}\to \mathbf{u}^*$ in $L^\infty(Q_{T^*};\mathbb{R}^4)$, and $M:K\to K$ is continuous. \hfill$\square$

\subsection{Proof of Lemma 2.3.2 part (b)}
We first recall the below Gagliardo-Nirenberg inequality.

\medskip
\noindent\textbf{Proposition 5.1.1} (Gagliardo–Nirenberg inequality on bounded domains).\\
Let $\Omega \subset \mathbb{R}^n$ be a measurable, bounded, open, and connected domain.  
Let $1 \le q \le \infty$ and $1 \le r \le \infty$.  
Let $j,m \in \mathbb{N}$ with $j<m$, and let $p \ge 1$.  
Assume $m-j-\tfrac{n}{r}$ is a nonnegative integer and let $\theta \in [0,1]$ such that the relations
$$
\frac{1}{p} \,=\, \frac{j}{n} \,+\, \theta \Big(\frac{1}{r}-\frac{m}{n}\Big) \,+\, \frac{1-\theta}{q},\quad \frac{j}{m} \le \theta < 1
$$
hold. Then there exists $C_g>0$, depending only on $n,m,j,p,q,r,\theta$ and $\Omega$, such that for all $u$ with $D^m u \in L^r(\Omega)$ and $u \in L^q(\Omega)$,
\begin{equation} \label{GNI}
\| D^j u\|_{L^p(\Omega)} 
\le C_g \,\| D^m u\|_{L^r(\Omega)}^\theta \,\| u\|_{L^q(\Omega)}^{\,1-\theta} 
+ C_g \,\| u\|_{L^{\sigma}(\Omega)},
\end{equation}
where $\sigma \in [1,\infty]$ is arbitrary \cite{Nirenberg1959}. 

\medskip
With Proposition 5.1.1 in hand, we move on to the proof.

\medskip
\noindent\textit{Proof (Lemma~2.3.2(b)).} Assume Lemma~2.3.2 \textbf{(a)} holds, i.e.,
\begin{equation}\label{A:L2assump}
\sup_{0\le t\le T^*}\int_\Omega u(x,t)^2\,dx \le M_0(T^*).
\end{equation}
Fix $p\in\mathbb{N}$, $p\ge 3$, and set
$$
w:=u^{\frac{p+1}{2}}\qquad(\text{so } w^2=u^{p+1}).
$$

\medskip
\noindent\textbf{Step 1: an energy for $w$.}
Multiply \eqref{pdes1} by $u^p$ and integrate over $\Omega$. Using \eqref{pdes2} and integration by parts,
\begin{align}\label{u_p}
\frac{1}{p+1}\frac{d}{dt}\int_\Omega u^{p+1}\,dx
+ \underbrace{p\int_\Omega D\nabla u\cdot \nabla u\,u^{p-1}\,dx}_{I}
&= \underbrace{\int_\Omega h\cdot \nabla u\,u^p\,dx}_{II}
+ \underbrace{\int_\Omega f(x,t,u)\,u^p\,dx}_{III}.
\end{align}
Rewrite $I$ in terms of $w$:
\begin{equation}\label{u_p1}
I \ge \frac{4p\,d_0}{(p+1)^2}\int_\Omega |\nabla w|^2\,dx .
\end{equation}
For $II$, since $\nabla(w^2)=(p+1)u^p\nabla u$, Young's inequality with $\epsilon_1>0$ gives
\begin{equation}\label{u_p2}
II
= \frac{1}{p+1}\int_\Omega h\cdot\nabla(w^2)\,dx
\le \frac{h^*}{p+1}\Big(\epsilon_1\!\int_\Omega|\nabla w|^2dx+\frac{1}{\epsilon_1}\!\int_\Omega w^2dx\Big).
\end{equation}
For $III$, use Young with $\epsilon_2>0$ (with $s=\frac{p+1}{p-1}$, $r=\frac{p+1}{2}$) by letting
$a=u^{p-1}$ and $b=fu$:
\begin{align}\label{L^2_V_1}
III
&= \int_{\Omega} f\,u\,u^{p-1}\,dx\\
&\le\,\frac{\epsilon_2(p-1)}{p+1}\int_{\Omega}u^{p+1}\,dx
+\frac{2}{p+1}\,\epsilon_2^{-\frac{p-1}{2}}\underbrace{\int_{\Omega}(fu)^{\frac{p+1}{2}}\,dx}_{III_1}.\nonumber
\end{align}
Using $fu\le \bar C(1+u^2)$, we obtain
\begin{equation}\label{L^2_V_2}
III_1 \le \bar{C}^{\frac{p+1}{2}}\int_{\Omega}u^{p+1}\,dx + \bar{C}^{\frac{p+1}{2}}|\Omega|.
\end{equation}
Hence
\begin{equation}\label{L^2_V_3}
III \le \Big(\frac{\epsilon_2(p-1)}{p+1}
+\frac{2}{p+1}\epsilon_2^{-\frac{p-1}{2}}\bar{C}^{\frac{p+1}{2}}\Big)\!\int_{\Omega}u^{p+1}\,dx
+\frac{2}{p+1}\epsilon_2^{-\frac{p-1}{2}}\bar{C}^{\frac{p+1}{2}}|\Omega|.
\end{equation}
Substitute \eqref{u_p1}, \eqref{u_p2}, and \eqref{L^2_V_3} into \eqref{u_p}, multiply by $(p+1)$, and choose
$\epsilon_1=\frac{2pd_0}{h^*(p+1)}$. Then
\begin{equation}\label{L^p_F}
\frac{d}{dt}\int_{\Omega}w^2\,dx + C_1(p)\int_{\Omega}|\nabla w|^2\,dx
\le C_2(p)\int_{\Omega}w^2\,dx + C_3(p),
\end{equation}
where $C_1(p)=\frac{2pd_0}{p+1}>0$ and $C_2(p),C_3(p)$ depend only on $p$, $\epsilon_2$, and known data.

\medskip
\noindent\textbf{Step 2: Gagliardo-Nirenberg on $w$.}
By Proposition~5.1.1 with $j=0$, $m=q=\sigma=1$, and $r=p=2$,
\begin{equation}\label{GNI_s}
\| w\|_{L^2(\Omega)}
\le C_g \|\nabla w\|_{L^2(\Omega)}^\theta \|w\|_{L^1(\Omega)}^{1-\theta}
+ C_g \|w\|_{L^1(\Omega)},
\end{equation}
where $\theta=\frac{n}{n+2}$. Using $(a+b)^2\le 2(a^2+b^2)$ and weighted Young's inequality (with $\epsilon_3>0$), we obtain
\begin{equation}\label{GNI_s2}
\|\nabla w\|_{L^2(\Omega)}^2
\ge C_4\|w\|_{L^2(\Omega)}^2 - C_5\|w\|_{L^1(\Omega)}^2,
\end{equation}
for constants $C_4,C_5>0$ depending only on $\Omega$ and $C_g$.

\medskip
\noindent\textbf{Step 3: Alikakos iteration.} (see \cite{Alikakos1979})
Inserting \eqref{GNI_s2} into \eqref{L^p_F} yields
\begin{equation}\label{it_base}
\frac{d}{dt}\int_{\Omega} w^2 \,dx
+\Big(\underbrace{C_1(p)C_4 - C_2(p)}_{C_6(p)}\Big)\int_{\Omega} w^2 \,dx
\le C_3(p) + C_1(p)C_5\|w\|_{L^1(\Omega)}^2.
\end{equation}

\medskip
\textit{Base step $p_1=3$.}\\
For $p=p_1=3$, we have $w=u^2$. From \eqref{A:L2assump}, we have
$$
\|w\|_{L^1(\Omega)}=\int_\Omega u^2\,dx \le M_0(T^*),
\quad\Rightarrow\quad
\|w\|_{L^1(\Omega)}^2 \le (M_0(T^*))^2.
$$
Choose $\epsilon_2,\epsilon_3$ so that $C_6(p_1=3)>0$, and apply Grönwall to \eqref{it_base} with $p=3$:
\begin{equation}\label{L4}
\sup_{0 \le t \le T^*}\int_\Omega w^{2}\,dx
= \sup_{0 \le t \le T^*}\int_\Omega u^{4}\,dx \le \hat{M}(T^*)_{p_1+1} = \hat{M}(T^*)_4,
\end{equation}
where $\hat{M}(T^*)_{4}$ depends only on $M_0(T^*)$ and known data.

\medskip
\textit{Iteration on $p_{k+1}$.}\\
Following Alikakos [\cite{Alikakos1979}], set $p_1=3$ and $p_{k+1}:=2p_k+1$ for $k=1,2,\dots$, so that $p_{k}=2^{k+1}-1\to\infty$. For each $k$, define the $(k+1)$-th iteration  $w_{k+1}:=u^{\frac{p_{k+1}+1}{2}}=u^{p_k+1}$, so that
\begin{equation}\label{k+1_bound}
\|w_{k+1}\|_{L^1(\Omega)}=\int_\Omega u^{\,p_k+1}\,dx= \int_\Omega u^{2^{k+1}}\,dx.
\end{equation}
Assume the induction hypothesis
$$
\sup_{0 \le t \le T^*}\int_\Omega u^{\,p_k+1}\,dx \le \hat{M}(T^*)_{p_k+1},
$$
so
$$
\sup_{0 \le t \le T^*}\|w_{k+1}\|_{L^1(\Omega)} \le \hat{M}(T^*)_{p_k+1}
\quad\Rightarrow\quad
\sup_{0 \le t \le T^*}\|w_{k+1}\|_{L^1(\Omega)}^{2} \le \big(\hat{M}(T^*)_{p_k+1}\big)^2.
$$
Using \eqref{it_base} with $p=p_{k+1}$,
\begin{equation}\label{it_p+1}
\frac{d}{dt}\int_\Omega w_{k+1}^2dx
+ C_6(p_{k+1})\int_\Omega w_{k+1}^2dx
\le C_3(p_{k+1}) + C_1(p_{k+1})C_5\|w_{k+1}\|_{L^1(\Omega)}^2.
\end{equation}
Choose $\epsilon_2,\epsilon_3$ so $C_6(p_{k+1})>0$; by Grönwall,
\begin{equation}\label{k+1_com_inc}
\sup_{0 \le t \le T^*}\int_\Omega w_{k+1}^2 dx
=\sup_{0 \le t \le T^*}\int_\Omega u^{2p_k+2} dx
\le \hat{M}(T^*)_{p_{k+1}+1}.
\end{equation}
Therefore, for all $k$,
$$
\sup_{0\le t\le T^*}\|u(\cdot,t)\|_{L^{\,p_k+1}(\Omega)} \le C_k(T^*),
$$
for a sequence $p_k\to\infty$. Since $\Omega$ is bounded, letting $k\to\infty$ gives
$$
\sup_{0\le t\le T^*}\|u(\cdot,t)\|_{L^\infty(\Omega)} \le \hat M(T^*).
$$
\hfill$\square$

\bibliography{biblio}

\end{document}